\documentclass[12pt]{article}
\usepackage{amsmath, amssymb,latexsym}
\usepackage{color}
\usepackage{tikz}
\usetikzlibrary{shapes.misc}
\tikzset{cross/.style={cross out, draw=black, minimum size=2*(#1-\pgflinewidth), inner sep=0pt, outer sep=0pt},
cross/.default={5pt}}

\title{Spectra for semiclassical operators with periodic bicharacteristics in dimension two}
\author{Michael A. Hall\\Department of Mathematics \\University of California \\ Los Angeles
\\CA 90095-1555, USA\\\small michaelhall@math.ucla.edu \and
Michael Hitrik\\Department of Mathematics \\University of California \\ Los Angeles
\\CA 90095-1555, USA\\\small hitrik@math.ucla.edu \and
Johannes Sj\"ostrand\\IMB, Universit\'e de Bourgogne\\9, Av. A. Savary, BP 47870\\FR--21078 Dijon, France \\and UMR 5584 CNRS
\\\small johannes.sjostrand@u-bourgogne.fr}
\date{}

\def\wrtext#1{\relax\ifmmode{\leavevmode\hbox{#1}}\else{#1}\fi}
\def\abs#1{\left|#1\right|}
\def\begeq{\begin{equation}}
\def\endeq{\end{equation}}

\def\neigh{neighborhood}

\def\Re{{\rm Re\,}}
\def\Im{{\rm Im\,}}

\newcommand{\eps}{\varepsilon}
\def\part#1{\frac{\partial}{\partial #1}}

\def\norm#1{||\,#1\,||}

\newcommand{\real}{\mbox{\bf R}}
\newcommand{\comp}{\mbox{\bf C}}
\newcommand{\z}{\mbox{\bf Z}}
\newcommand{\nat}{\mbox{\bf N}}

\renewcommand{\Re}{\mbox{\rm Re\,}}
\renewcommand{\Im}{\mbox{\rm Im\,}}

\renewcommand{\exp}{\mbox{\rm exp\,}}

\newtheorem{dref}{Definition}[section]

\newtheorem{theo}[dref]{Theorem}
\newtheorem{prop}[dref]{Proposition}

\begin{document}
\maketitle

\begin{abstract}
We study the distribution of eigenvalues for selfadjoint $h$--pseudodifferential operators in dimension two, arising as perturbations of selfadjoint operators
with a periodic classical flow. When the strength $\eps$ of the perturbation is $\ll h$, the spectrum displays a cluster structure, and assuming that $\eps \gg h^2$
(or sometimes $\gg h^{N_0}$, for $N_0 >1$ large), we obtain a complete asymptotic description of the individual eigenvalues inside subclusters, corresponding to the
regular values of the leading symbol of the perturbation, averaged along the flow.
\end{abstract}
\tableofcontents

\section{Introduction and statement of results}
\setcounter{equation}{0}
The spectral theory of selfadjoint partial differential operators, whose associated classical flow is periodic, has a long and distinguished tradition,
starting with the classical works of J. J. Duistermaat-V. Guillemin~\cite{DG} and A. Weinstein~\cite{Wei}, in the high energy limit,
in the case of compact manifolds. Subsequently, many important contributions to the theory were given,~\cite{CdV},~\cite{Gu},~\cite{BdMG},~\cite{Zel},
and the case of semiclassical pseudodifferential operators was treated in~\cite{HeRo},~\cite{Dozias},~\cite{Ivrii}. In particular, assuming that the Hamilton flow
is periodic in some energy shell, the cluster structure of the spectrum has been established in~\cite{HeRo}. That work also contains some precise results
concerning the semiclassical asymptotics for the counting function of eigenvalues in the clusters, with the celebrated Bohr-Sommerfeld quantization rule
obtained as a special case in dimension one, see also~\cite{DS}. Let us also remark that apart from their intrinsic interest in spectral theory, starting from
the work~\cite{DG}, operators with periodic classical flow have frequently served as a source of examples of situations where various spectral estimates become
optimal --- see~\cite{BSSY} for a recent manifestation of this in the context of uniform $L^p$ resolvent estimates for the Laplacian on a compact Zoll manifold.

\bigskip
\noindent
The purpose of this paper is to show how the microlocal techniques of \cite{HS1}, \cite{HS2}, \cite{HS3}, developed in the context of analytic
non-selfadjoint perturbations of selfadjoint operators with periodic classical flow in dimension two, apply to a class of {\it selfadjoint} operators of the form
$P_{\eps} = P(x,hD_x) + \eps Q(x,hD_x)$. Here $P=P(x,hD_x)$ and $Q=Q(x,hD_x)$ are selfadjoint, with $P$ elliptic at infinity and with the classical flow of $P$
periodic in a band of energies around $0$. It is then well-known, and will be recalled below, that the spectrum of $P$ near $0$ exhibits a cluster structure,
each cluster being of size $\leq {\cal O}(h^2)$, and with the separation between the adjacent clusters of size $h$. The parameter $\eps\in [0,\eps_0)$, $\eps_0>0$,
measures the strength of the selfadjoint perturbation, and in order to have the clustering for the spectrum of $P_{\eps}$, one should have $\eps \ll h$. The general
problem is then to understand the internal structure of the spectral clusters of the perturbed operator $P_{\eps}$ in some detail, in the semiclassical limit
$h\rightarrow 0$. In this work, assuming that $\eps \gg h^{N_0}$, where $N_0 \geq 2$ depends on the size of the spectral clusters of $P$, we obtain semiclassical complete asymptotic
expansions for the individual eigenvalues of $P_{\eps}$ in subclusters, corresponding to regular values of the leading symbol of the perturbation $Q$, averaged along the classical flow. We remark that
contrary to~\cite{HS1},~\cite{HS2},~\cite{HS3}, no analyticity assumptions are needed here, and the spectral analysis is carried out within the framework of the
standard $L^2$--spaces.

\bigskip
\noindent
Let $M$ stand for $\real^2$ or a smooth compact $2$--dimensional Riemannian manifold. When $M = \real^2$, we let
\begeq
\label{eq1.1}
P_{\eps} = P^w(x,hD_x,\eps;h),\quad 0 < h \ll 1,
\endeq
be the $h$--Weyl quantization on $\real^2$ of a symbol $P(x,\xi,\eps;h)$ depending smoothly on $\eps \in {\rm neigh}(0,\real)$ taking values in the
symbol class $S(m)$. Here $m$ is assumed to be an order function on $\real^{4}$ in the sense that $m>0$ and
\begeq
\label{eq1.2}
\exists C_0 \geq 1, N_0>0,\ m(X)\leq C_0 \langle{X-Y\rangle}^{N_0} m(Y),\ X,Y\in \real^{4}.
\endeq
The symbol class $S(m)$ is given by
$$
S(m)=\left\{ a\in C^{\infty}(\real^{4}); \forall \alpha \in \nat^{4}, \exists C_{\alpha}>0, \forall X \in \real^{4},
\  |\partial_X^{\alpha} a(X)| \leq C_{\alpha} m(X)\right\}.
$$
We shall assume throughout that
\begeq
\label{eq1.3}
m\geq 1.
\endeq

\medskip
\noindent
Assume furthermore that
\begeq
\label{eq1.4}
P(x,\xi,\eps;h) \sim \sum_{j=0}^{\infty} h^j p_j(x,\xi,\eps), \quad h\to 0,
\end{equation}
in the space $S(m)$. We make the ellipticity assumption,
\begeq
\label{eq1.5}
\abs{p_0(x,\xi,\eps)} \geq \frac{1}{C} m(x,\xi),\quad \abs{(x,\xi)} \geq C,
\endeq
for some $C>0$.

\bigskip
\noindent
When $M$ is a compact manifold, we first recall the standard class of semiclassical symbols on $T^*M$,
$$
S^m(T^*M) = \left\{ a(x,\xi;h) \in C^{\infty}(T^*M \times (0,1]): \abs{\partial_x^{\alpha}\partial_{\xi}^{\beta} a(x,\xi;h)} \leq C_{\alpha\beta}
\langle{\xi\rangle}^{m-\abs{\beta}}\right \}.
$$
Associated to $S^m(T^*M)$ is the corresponding class of semiclassical pseudodifferential operators denoted by $L^m(M)$. Let $P_{\eps}$ be a $C^{\infty}$ function
of $\eps \in {\rm neigh}(0,\real)$ with values in $L^m(M)$, $m>0$. Let $\widetilde{M} \subset M$ be a coordinate chart identified with a convex
bounded domain in $\real^n$ in such a way that the Riemannian volume element $\mu(dx)$ reduces to the Lebesgue measure in $\widetilde{M}$. We then have on
$\widetilde{M}$, for every $u\in C^{\infty}_0(\widetilde{M})$,
\begeq
\label{eq1.6}
P_{\eps} u(x) = \frac{1}{(2 \pi h)^2} \int\!\!\!\int e^{\frac{i}{h}(x-y)\cdot \xi} p\left(\frac{x+y}{2},\xi,\eps;h\right) u(y)dy\,d\xi + Ru(x).
\endeq
Here $p(x,\xi,\eps;h)$ is a smooth function of $\eps$ with values in $S^m_{{\rm loc}}(\widetilde{M}\times \real^2)$, and $R$ is negligible in the sense
that its Schwartz kernel $R(x,y)$ satisfies $\partial_x^{\alpha} \partial_y^{\beta} R(x,y) = {\cal O}(h^{\infty})$, for all $\alpha$, $\beta$. We further assume
that the symbol $p(x,\xi,\eps;h)$ has an asymptotic expansion in $S^m_{{\rm loc}}(\widetilde{M}\times \real^2)$, as $h\rightarrow 0$,
\begeq
\label{eq1.61}
p(x,\xi,\eps;h) \sim \sum_{j=0}^{\infty} h^j p_j(x,\xi,\eps),\quad p_j \in S^{m-j}_{{\rm loc}}.
\endeq
The semiclassical principal symbol of $P_{\eps}$ in this case is given by $p_0(x,\xi,\eps)$, and we make the ellipticity assumption,
$$
\abs{p_0(x,\xi,\eps)} \geq \frac{1}{C} \langle{\xi\rangle}^m,\quad (x,\xi)\in T^*M,\quad \abs{\xi}\geq C,
$$
for some $C>0$. Here we recall that since $M$ has been equipped with some Riemannian metric, $\abs{\xi}$ and $\langle{\xi\rangle} = (1+\abs{\xi}^2)^{1/2}$ are
well defined. Let us also recall from~\cite{SjZw02} that while the complete symbol $p$ in (\ref{eq1.61}) depends on the choice of local coordinates, the principal symbol
$p_0$ and the subprincipal symbol $p_1$ are invariantly defined, provided that we use local coordinates in (\ref{eq1.6}) for which the Riemannian volume density
becomes equal to the Lebesgue measure.

\medskip
\noindent
In what follows, we shall write $p_{\eps}$ for the principal symbol $p_0(x,\xi,\eps)$ of $P_{\eps}$, and simply $p$ for $p_0(x,\xi,\eps=0)$. We shall assume
that for $\eps \in {\rm neigh}(0,\real)$,
\begeq
\label{eq1.7}
P_{\eps}\,\,\, \wrtext{is formally selfadjoint}.
\endeq
In the case when $M$ is compact, we let the underlying Hilbert space be $L^2(M,\mu(dx))$.

\bigskip
\noindent
For $h>0$ small enough and when equipped with the domain $H(m)$, the naturally defined Sobolev space associated with the order function $m$
(so that in the compact case, $H(m)$ is the standard semiclassical Sobolev space $H^m(M)$), $P_{\eps}$ becomes a well-defined selfadjoint operator on $L^2(M)$.
Moreover, the assumptions above imply that the spectrum of $P_{\eps}$ in a fixed neighborhood of $0$ is discrete, for $h>0$ and $\eps\geq 0$ small enough.

\medskip
\noindent
We shall assume that the energy surface $p^{-1}(0)\subset T^*M$ is connected and that $dp\neq 0$ along $p^{-1}(0)$. Let
$H_p = p'_{\xi}\cdot \frac{\partial}{\partial x} - p'_x\cdot \frac{\partial}{\partial \xi}$ be the Hamilton vector field of $p$. We introduce the following
basic assumption, assumed to hold throughout this work: for $E\in {\rm neigh}(0,\real)$,
\begin{eqnarray}
\label{eq1.8}
\hbox{The }H_p\hbox{-flow is periodic on }p^{-1}(E) \hbox{ with minimal}\\
\hbox{period }T(E)>0 \hbox{ depending smoothly on } E.\nonumber
\end{eqnarray}

\medskip
\noindent
Let $g: {\rm neigh}(0,\real) \to \real$ be the smooth function defined by
\begeq
\label{eq1.85}
g'(E) = \frac{T(E)}{2\pi},\quad g(0) = 0,
\endeq
so that $g\circ p$ has a $2\pi$-periodic Hamilton flow. Set $f = g^{-1}$. We then have the following well known result, due to~\cite{HeRo}, following the earlier
works~\cite{Wei},~\cite{CdV}.

\begin{theo}
Assume that the subprincipal symbol of $P_{0}$ vanishes. Then the spectrum of $P_{0}$ near $0$ is contained in the union of the intervals of the form
\begeq
\label{eq1.86}
I_k = f(h(k-\theta)) + [-{\cal O}(h^2), {\cal O}(h^2)],\quad k\in \z,
\endeq
pairwise disjoint for $h>0$ small enough. Here $\theta = \alpha_1/4 + S_1/2\pi h$, where $\alpha_1 \in \z$ and $S_1\in \real$ are the Maslov index and the classical
action, respectively, computed along a closed $H_p$-trajectory $\subset p^{-1}(0)$, of period $T(0)$.
\end{theo}

\medskip
\noindent
{\it Remark}. We refer to~\cite{M} for a self-contained discussion of Maslov indices of loops of Lagrangian subspaces and closed Hamiltonian trajectories.

\medskip
\noindent
{\it Remark}. Let us observe that up to a constant, the function $g$ in (\ref{eq1.85}) is equal to $1/2\pi$ times the classical action along a closed $H_p$--orbit
in $p^{-1}(E)$, $E\in {\rm neigh}(0,\real)$, see~\cite{DS}.

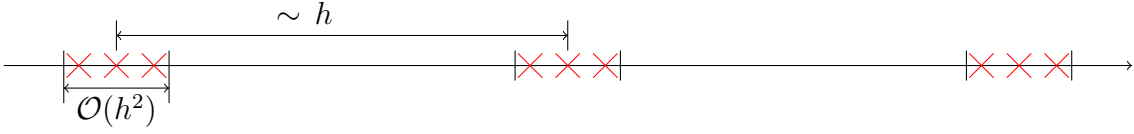
\begin{figure} [ht]
\centering
\begin{tikzpicture}
\draw [->] (-1,0) -- (14,0);
\draw (0,0) node[cross,red] {};
\draw (0.5,0) node[cross,red] {};
\draw (1,0) node[cross,red] {};
\draw [-] (-0.2,-0.5) -- (-0.2,0.2);
\draw [-] (1.2,-0.5) -- (1.2,0.2);
\draw [<->]  (-0.2,-0.3) -- (1.2,-0.3);
\node at (0.5,-0.6) {$\mathcal{O}(h^2)$};

\draw (6,0) node[cross,red] {};
\draw (6.5,0) node[cross,red] {};
\draw (7,0) node[cross,red] {};
\draw [-] (5.8,-0.2) -- (5.8,0.2);
\draw [-] (7.2,-0.2) -- (7.2,0.2);

\draw (12,0) node[cross,red] {};
\draw (12.5,0) node[cross,red] {};
\draw (13,0) node[cross,red] {};
\draw [-] (11.8,-0.2) -- (11.8,0.2);
\draw [-] (13.2,-0.2) -- (13.2,0.2);

\draw [-] (0.5,0.2) -- (0.5,0.6);
\draw [-] (6.5,0.2) -- (6.5,0.6);
\draw [<->] (0.5,0.4) -- (6.5,0.4);
\node at (3,0.7) {$\sim\, h$};
\end{tikzpicture}
\caption{Spectral clusters for the unperturbed operator $P_0$. The size of each cluster is ${\cal O}(h^2)$, with the separation between adjacent clusters being
of order $h$.}
\label{fig_1}
\end{figure}

\bigskip
\noindent
Let us write
\begeq
\label{eq1.9}
p_{\eps} = p + \eps q + {\cal O}(\eps^2 m),
\endeq
in the case $M = \real^2$, and $p_{\eps} = p + \eps q + {\cal O}(\eps^2 \langle{\xi\rangle}^m)$ in the compact case. Here $q$ is smooth and
real-valued on $T^*M$. Let
\begeq
\label{eq1.10}
\langle{q\rangle} = \frac{1}{T(E)}\int_{-T(E)/2}^{T(E)/2} q \circ \exp(tH_p)\,dt \quad \wrtext{on}\,\, p^{-1}(E),
\endeq
and notice that $H_p \langle{q\rangle} = 0$. In view of (\ref{eq1.8}), the space of closed $H_p$--orbits in $p^{-1}(0)$,
$$
\Sigma = p^{-1}(0)/\exp(\real H_p),
$$
is a 2-dimensional symplectic manifold, and $\langle{q\rangle}$ can naturally be viewed as a function on $\Sigma$.

\medskip
\noindent
The following is the main result of this work.

\begin{theo}
Let us assume that {\rm (\ref{eq1.8})} holds and that the subprincipal symbol of $P_{0}$ vanishes. Let $F_0$ be a regular value of $\langle{q}\rangle$, considered
as a function on $\Sigma$. Assume that $\langle{q\rangle}^{-1}(F_0)\subset \Sigma$ is a connected closed curve, and let us introduce the corresponding Lagrangian
torus,
$$
\Lambda_{0,F_0} = \{\rho\in T^*M;\,\, p(\rho)  = 0,\,\,\, \langle{q\rangle}(\rho)= F_0\}.
$$
When $\gamma_1$ and $\gamma_2$ are the fundamental cycles in $\Lambda_{0,F_0}$ with $\gamma_1$ being given by a closed $H_p$--trajectory of minimal period, we
write $S=(S_1,S_2)$ and $\alpha=(\alpha_1,\alpha_2)$ for the actions and the Maslov indices of the cycles, respectively. Assume next that the spectrum of $P_0$
near $0$ clusters into bands of size ${\cal O}(1)h^{N_0}$, for some $N_0 \geq 2$. Let us assume that
$$
h^{N_0} \ll \eps \ll h.
$$
Let $C>0$ be large enough. There exists a smooth function
\begeq
\label{eq1.10.1}
f(\xi_1;h) = f(\xi_1) + \sum_{j=2}^{N_0-1} h^j f_j(\xi_1), \quad \xi_1\in {\rm neigh}(0,\real),
\endeq
such that for each $k \in \z$, with $h(k-\alpha_1/4) - S_1/2\pi$ small enough, the eigenvalues of $P_{\eps}$ in the set
\begeq
\label{eq1.11}
\abs{z-f\left(h(k-\frac{\alpha_1}{4})-\frac{S_1}{2\pi};h\right)-\eps F_0} < \frac{\eps}{C}
\endeq
are given by
\begeq
\label{eq1.12}
\widehat{P}\left(h(k-\frac{\alpha_1}{4})-\frac{S_1}{2\pi},h(\ell-\frac{\alpha_2}{4})-\frac{S_2}{2\pi},\eps,\frac{h^{N_0}}{\eps};h\right)
+{\cal O}(h^{\infty}),\quad \ell \in \z.
\endeq
Here $\widehat{P}(\xi,\eps,h^{N_0}/\eps;h)$ is smooth in $\xi \in {\rm neigh}(0,\real^2)$, smooth in
$\eps,\frac{h^{N_0}}{\eps}\in {\rm neigh}(0,\real)$, and has a complete asymptotic expansion in the space of such functions, as $h\rightarrow 0$,
$$
\widehat{P}\left(\xi,\eps,\frac{h^{N_0}}{\eps};h\right)\sim f(\xi_1;h)+\eps \left(r_0\left(\xi,\eps,\frac{h^{N_0}}{\eps}\right)+h
r_1\left(\xi,\eps,\frac{h^{N_0}}{\eps}\right)+\ldots\right).
$$
We have
$$
r_0(\xi)=\langle{q}\rangle(\xi)+{\cal O}\left(\eps+\frac{h^{N_0}}{\eps}\right),\quad r_j={\cal O}(1),\quad j\geq 1,
$$
corresponding to the action-angle coordinates near the Lagrangian torus $\Lambda_{0,F_0}$.
\end{theo}

\begin{figure} [ht]
\centering
\begin{tikzpicture}
\draw [ ->] (-1,0) -- (14,0);
\draw [ -] (-0.2,-0.9) -- (-0.2,0.2);
\draw [ -] (13.2,-0.9) -- (13.2,0.2);
\draw [<->]  (-0.2,-0.8) -- (13.2,-0.8);
\node at (6,-1.1) {$\mathcal{O}(\varepsilon)$};
\draw [-] (6.5,-0.2) -- (6.5,0.2);
\node at (6.2,0.5) {$f(h(k-\theta))$};

\draw (9,0) node[cross,red] {};
\draw (9.5,0) node[cross,red] {};
\draw (10,0) node[cross,red] {};
\draw (10.5,0) node[cross,red] {};
\draw [ -] (10.5,-0.2) -- (10.5,0.2);
\node at (10.5,0.5) {$f(h(k-\theta))+\eps F_0$};
\draw (11,0) node[cross,red] {};
\draw (11.5,0) node[cross,red] {};
\draw (12,0) node[cross,red] {};

\draw [ -] (9.5,-0.1) -- (9.5,-0.5);
\draw [ -] (10,-0.1) -- (10,-0.5);
\draw [<->]  (9.5,-0.3) -- (10,-0.3);
\node at (10,-0.6) {$\sim\varepsilon h$};

\node at (8.8,0) {$($};
\node at (12.2,0) {$)$};
\end{tikzpicture}
\caption{Spectral asymptotics in a subcluster of the $k$th spectral cluster of the operator $P_{\eps}$, corresponding to the regular value $F_0$ of the leading
symbol of the perturbation, averaged along the classical flow. Here $h^2 \ll \eps \ll h$. The red crosses represent the eigenvalues of $P_{\eps}$ in (\ref{eq1.11}),
given by (\ref{eq1.12}).}
\label{fig_2}
\end{figure}
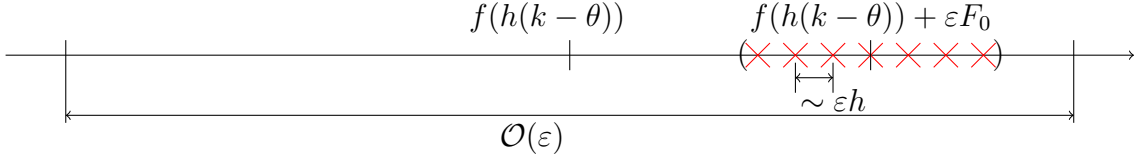

\medskip
\noindent
{\it Remark}. In the case when the compact manifold $\Lambda_{0,F_0}$ has several connected components, the result of Theorem 1.2 may be extended by showing that
the set of eigenvalues $z$ in (\ref{eq1.11}) agrees with the union of the spectral contributions coming from each of the connected components,
modulo ${\cal O}(h^{\infty})$, each contribution having the description as in the theorem. See also the discussion in Section 3 below.

\medskip
\noindent
{\it Example}. Let $M$ be a compact symmetric surface of rank one, and let $P_{0} = -h^2 \Delta - 1$ on $M$. The assumption (\ref{eq1.8}) then holds and
from~\cite{Gu} we know that the spectrum of $P_0$ clusters into bands of diameter $0$ and separation of order $h$. Furthermore, the eigenvalues of $P_0+1$ depend
quadratically on $h$, and we may conclude that the functions $f_j$ in (\ref{eq1.10.1}) vanish, for $j>2$, while $f_2$ is a constant. Taking
$2 < N_0 \in \nat$ to be any fixed integer, we see that the result of Theorem 1.2 applies to
the Schr\"odinger operator $P_{\eps} = P_{0} + \eps q = h^2(-\Delta +q)-1$, where $\eps = h^2$ and $q \in C^{\infty}(M)$ is
real-valued, cf. with~\cite{Wei},~\cite{CdV1},~\cite{Gr}. Let now $A$ be a smooth real-valued $1$-form on $M$ and consider the magnetic Schr\"odinger operator on
$M$ given by
\begeq
\label{eq1.13}
P_{\eps} = h^2 d_{\lambda A}^* d_{\lambda A} + h^2 q-1.
\endeq
Here $d_A u  = du + iA\wedge u$, and $d_A^*$ is the Riemannian adjoint of $d_A$. Theorem 1.2 applies to the operator $P_{\eps}$ in
(\ref{eq1.13}) when $\eps = \lambda h$, provided that $0 < \lambda \ll 1$ is sufficiently small but fixed. In this case,
$q(x,\xi) = q_A(x,\xi) = 2\langle{\xi,A^{\sharp}\rangle}$, where $A^{\sharp}$ is the vector field associated to $A$ by means of the Riemannian metric. We may
also remark if $B=dA$ is the magnetic field and $dA = d\widetilde{A}$, then, since $H^1(M) = 0$, we have $\widetilde{A} = A + d\varphi$, where $\varphi\in C^{\infty}(M)$ is
real-valued. It follows that $q_{\widetilde{A}} - q_A = H_p\varphi$, where $p = \xi^2$ is the leading symbol of $P_0$, and therefore
the flow average $\langle{q_{\widetilde A}\rangle} = \langle{q_A\rangle}$ depends on the magnetic 2-form $B = dA$ only. See also Section 6 below.

\medskip
\noindent
{\it Remark}. It seems quite likely that the result of Theorem 1.2 can be extended to the case when $F_0$ is a non-degenerate critical value of
$\langle{q\rangle}$, cf. with~\cite{CdV1},~\cite{Gr}. We would also like to mention that the result of Theorem 1.2 can be viewed as a Bohr-Sommerfeld quantization
condition in the spectral clusters, corresponding to the regular values for a reduced one-dimensional operator, and here there are some direct links
with~\cite{ChVu},~\cite{Ch}, and the theory of Toeplitz operators on reduced compact symplectic spaces such as $\Sigma$. See also~\cite{GUW}.

\medskip
\noindent
The plan of the paper is as follows. In Section 2, after re-deriving the clustering of the spectrum of $P_0$, we carry out an averaging reduction of
$P_{\eps}$, microlocally in an energy shell. In Section 3, we microlocalize further to a suitable Lagrangian torus and construct a quantum Birkhoff normal form for
$P_{\eps}$ near the torus, very much following the approach of~\cite{HS1}. In Section 4 we solve a suitable global Grushin problem for $P_{\eps}$
and identify the spectrum in the subclusters precisely, thereby completing the main part of the proof of Theorem 1.2. In Section 5, we complete the discussion by
addressing the case when the spectral clusters of $P_0$ are of size ${\cal O}(h^{N_0})$, $N>2$, and show how to reach smaller values of $\eps$ in Theorem 1.2 in
this case. In Section 6, we finally give an application to the magnetic Schr\"odinger operator on $\real^2$ in the resonant case.

\medskip
\noindent
{\bf Acknowledgements}. The first author was supported by the UCLA Dissertation Year Fellowship. The third author was supported by the grant NOSEVOL ANR 2011 BS 01019 01.

\section{Clustering of eigenvalues and averaging reduction}
\label{clustering}
\setcounter{equation}{0}
For future reference, it will be convenient and natural for us to start by recalling a proof of Theorem 1.1 --- see also Proposition 2.1 of~\cite{HS2}.
When $z\in {\rm neigh}(0,\real)$, let us consider the equation
\begeq
\label{eq2.1}
(P_0 -z)u = v,\quad u\in H(m).
\endeq

\medskip
\noindent
Let $\chi \in C_0^\infty(T^*M;[0,1])$ be such that $\chi = 1$ near $p^{-1}(0)$. Semiclassical elliptic regularity gives, with the $L^2$ norms throughout, that
\begeq
\label{eq2.2}
\norm{(1-\chi)u}  \leq \mathcal{O}(1)\norm{v} + \mathcal{O}(h^\infty)\norm{u},
\endeq
where $\chi = {\rm Op}^w_h(\chi)$ is the corresponding quantization. Here and in what follows, when $M = \real^2$, we use the $h$--Weyl quantization,
while when $M$ is compact, we fix the choice of the quantization map ${\rm Op}^w_h: S^m(T^*M) \rightarrow L^m(M)$, given by the Weyl quantization in special local
coordinates as recalled in the introduction, with the associated symbol map: $L^m(M) \rightarrow S^m(T^*M)/h^2 S^{m-2}(T^*M)$.

\bigskip
\noindent
Turning the attention to a neighborhood of $p^{-1}(0)$, let $\gamma\subset p^{-1}(0)$ be a closed $H_p$-orbit, where we know that $T(0)$ is the minimal period of
$\gamma$. From Section 3 of~\cite{HS1}, we recall the following result.

\begin{prop}
There exists a smooth real-valued canonical transformation
\begeq
\kappa:\, {\rm neigh}(\gamma,T^*M) \to {\rm neigh}(\tau = x = \xi = 0,T^*(S^1_t\times \real_x)),
\endeq
mapping $\gamma$ onto $\{\tau = x = \xi=0\}$, such that $p\circ \kappa^{-1} = f(\tau)$. Here $f$ has been defined in {\rm (\ref{eq1.85})}.
\end{prop}

\medskip
\noindent
Following~\cite{HS1}, we recall that the canonical transformation $\kappa$ can be implemented by a multi-valued microlocally unitary $h$-Fourier
integral operator $U = {\cal O}(1): L^2(M)\to L^2_f(S^1\times \real)$, so that the improved Egorov property holds --- see the discussion in Section 2
of \cite{HS1}. Here $L^2_f(S^1\times \real)$ is the space of functions defined microlocally near $\tau = x = \xi =0$ in $T^*(S^1 \times \real)$, which
satisfy the Floquet-Bloch periodicity condition,
\begeq
\label{eq2.3}
u(t-2\pi,x) = e^{2\pi i\theta}u(t,x),\quad \theta = \frac{S_1}{2\pi h} + \frac{\alpha_1}{4}.
\endeq
As explained in~\cite{HS1}, the multi-valuedness of $U$ is a reflection of the fact that the domain of definition of the canonical transformation $\kappa$ is not
simply connected, the corresponding first homotopy group being generated by the closed trajectory $\gamma$.

\medskip
\noindent
It follows that there exists a selfadjoint operator $\widetilde P$ with the leading symbol $f(\tau)$ near $\{\tau = x = \xi = 0\}$ and with vanishing
subprincipal symbol, so that $\widetilde{P} U = UP_0$ microlocally near $\gamma$, so that
\begeq
\label{eq2.4}
(\widetilde{P} U - UP_0){\rm Op}^w_h(\chi_1) = {\cal O}(h^{\infty}): L^2 (M)\rightarrow L^2(M),
\endeq
and
$$
\chi_2^w(x,hD_x) (\widetilde{P} U - UP_0) = {\cal O}(h^{\infty}): L^2_f (S^1\times \real)\rightarrow L^2_f (S^1\times \real),
$$
for every $\chi_1\in C^{\infty}_0({\rm neigh}(\gamma,T^*M))$ and for every $\chi_2\in C^{\infty}_0(T^*(S^1\times \real))$ supported near
$\tau=x=\xi=0$. The operator $\widetilde{P}$ acts on the space of functions satisfying the Floquet-Bloch condition (\ref{eq2.3}),
defined microlocally near $\tau = x = \xi=0$ in $T^*(S^1 \times \real)$.

\medskip
\noindent
Let us remark next that an orthonormal basis for the space $L^2_f(S^1)$ of functions $u\in L^2_{{\rm loc}}(\real)$ satisfying a Floquet-Bloch periodicity condition
analogous to (\ref{eq2.3}), with the $x$-variable suppressed, consists of the functions
\begeq
\label{eq2.5}
e_k(t) = \frac{1}{\sqrt{2\pi}} \exp(i(k-\theta)t),\quad \theta = \frac{S_1}{2\pi h} + \frac{\alpha_1}{4},\quad k \in \z,
\endeq
which satisfy $f(hD_t) e_k(t) = f(h(k-\theta)) e_k(t)$. It follows that if $z\in {\rm neigh}(0,\real)$ is such that
$$
\abs{z - f(h(k-\theta))} \geq C h^2,\quad k\in \z,
$$
for $C>1$ sufficiently large but fixed, then the operator
$$
\widetilde{P} - z =  f(hD_t) + h^2 R -z, \quad R = {\cal O}(1): L^2_f(S^1\times\real) \to L^2_f(S^1\times \real),
$$
is invertible, microlocally near $\tau=x=\xi=0$, with the norm of the inverse being ${\cal O}(h^{-2})$.

\medskip
\noindent
Let us now take finitely many closed $H_p$--trajectories $\gamma_1, \ldots ,\gamma_N\subset p^{-1}(0)$ and small open neighborhoods $\Omega_j$ of $\gamma_j$, with
$\Omega_j$ invariant under the $H_p$-flow, $1\leq j \leq N$, such that $p^{-1}(0) \subseteq \cup \Omega_j$. Associated to this open cover, we take cutoff
functions $0\leq \chi_j \in C_0^\infty(\Omega_j)$ such that $H_p \chi_j = 0$ and $\sum \chi_j  = 1$ near $p^{-1}(0)$. Let $U_j$ denote a multi-valued microlocally
unitary $h$-Fourier integral operator associated to the canonical transformation near $\gamma_j$, as in Proposition 2.1.

\medskip
\noindent
For each $j$, $1\leq j \leq N$, using (\ref{eq2.1}), we see that
\begeq
(P_0-z)\chi_j u + [\chi_j,P_0]u = \chi_j v,
\endeq
and therefore,
\begeq
\label{eq2.51}
U_j (P_0-z)\chi_j u = (\widetilde P-z)U_j \chi_j u + \mathcal{O}(h^\infty) u = U_j \chi_j v + U_j [P_0,\chi_j]u.
\endeq
When $z\in {\rm neigh}(0,\real)$ avoids the intervals $I_k$ in (\ref{eq1.86}), we just saw that the operator
$$
\widetilde{P} -z = f(hD_t) + h^2 R_j-z,
$$
possesses a microlocal inverse of norm $\mathcal{O}(1/h^2)$. Using (\ref{eq2.51}), we conclude that
\begeq
\label{eq2.52}
\norm{\chi_j u} \leq {\cal O}\left(\frac{1}{h^2}\right)\norm{v} + {\cal O}\left(\frac{1}{h^2}\right) \norm{[P,\chi_j]u} + {\cal O}(h^{\infty})\norm{u}.
\endeq
Since the subprincipal symbols of $P_0$ and $\chi_j$ both vanish, we have in the operator sense, $[P_0,\chi_j] = \mathcal{O}(h^3)$ --- see also~\cite{Gu0} for
composition rules for the subprincipal symbols. Using (\ref{eq2.52}) and summing over $j$ we get
\begeq
\label{eq2.6}
\norm{\sum_{j=1}^N \chi_j u} \leq {\cal O}\left(\frac{1}{h^2}\right) \norm{v} + \mathcal{O}(h)\norm{u}.
\endeq
Combining (\ref{eq2.2}) and (\ref{eq2.6}), we obtain
\begeq
\norm{u} \leq {\cal O}\left(\frac{1}{h^2}\right) \norm{v} + {\cal O}(h)\norm{u}.
\endeq
Taking $h$ small enough, we conclude that $(P_0-z)^{-1}$ exists and satisfies $(P_0-z)^{-1} = \mathcal{O}(1/h^2): L^2 \rightarrow L^2$, when
$z\in {\rm neigh}(0,\real)$ avoids the intervals $I_k$, $k\in \z$.

\bigskip
\noindent
Having recalled a proof of Theorem 1.1, let us now proceed to carry out an averaging reduction for the selfadjoint operator $P_{\eps}$, replacing the leading
symbol $q$ of the perturbation by its average along closed orbits of the $H_p$--flow. We notice that such a reduction has a very long
tradition~\cite{Wei},~\cite{Gu1}, and the following discussion will be therefore somewhat brief.

\bigskip
\noindent
Let $G_0 \in C^{\infty}$ in a neighborhood of $p^{-1}(0)$ be real-valued and such that
\begeq
\label{eq2.61}
H_p G_0 = q - \left\langle q \right\rangle,
\endeq
where $\langle{q\rangle}$ is the flow average, defined in (\ref{eq1.10}). As recalled in~\cite{HS1}, we may take
\begin{multline*}
G_0 \\ = \frac{1}{T(E)} \int_{-T(E)/2}^{T(E)/2}
\left[1_{\mathbf{R}_-}(t)\left(t + \frac{1}{2}T(E) \right) + 1_{\mathbf{R}_{+}}(t)\left(t - \frac{1}{2}T(E) \right)  \right]
q\circ\exp(tH_p) \, dt,
\end{multline*}
on $p^{-1}(E)$. By a Taylor expansion, we then get
$$
p_\eps \circ \exp(\eps H_{G_0}) = p + \eps (q - H_p G_0) + \mathcal{O}(\eps^2) = p + \eps \langle{q\rangle} + {\cal O}(\eps^2).
$$

\medskip
\noindent
Similarly, with $G_1,G_2, \ldots $ denoting a sequence of smooth real-valued functions to be determined, and $G \sim \sum_{j=0}^\infty \eps^j G_j$, if we
expand $p_\eps \circ \exp (\eps H_G)$ asymptotically, we claim that we can iteratively solve for $G_j$ so that
$$
p_\eps \circ \exp(\eps H_G) = p + \eps \left\langle q \right\rangle + \mathcal{O}(\eps^2),
$$
where the $\mathcal{O}(\eps^2)$ error term is real-valued and Poisson commutes with $p$, modulo $\mathcal{O}(\eps^\infty)$.
Explicitly, if $G_{\leq N} = G_0 + \eps G_1 + \eps^2 G_2 +  \ldots + \eps^{N} G_{N}$ satisfies
$$
p_\eps \circ \exp(\eps H_{G_{\leq N}}) =
p + \eps \left\langle q \right\rangle + \eps^2 q_2 + \ldots + \eps^{N+1} q_{N+1} + \eps^{N+2}r_{N+2} + \mathcal{O}(\eps^{N+3}),
$$
where $q_j$ are real-valued and $H_p q_j  = 0$, $2\leq j \leq N+1$, then for $G_{\leq N+1} = G_{\leq N} + \eps^{N+1}G_{N+1}$, with
$G_{N+1}\in C^\infty$ to be determined, we have by a variation on the Campbell-Hausdorff formula~\cite{Horm67},
\begin{multline}
\exp(\eps H_{G_{\leq N+1}}) = \exp(\eps H_{G_{\leq N}} + \eps^{N+2} H_{G_{N+1}}) \\
= \exp(\eps H_{G_{\leq N}})\exp(\eps^{N+2}H_{G_{N+1}})(1 + {\cal O}(\eps^{N+3})),
\end{multline}
where the $\mathcal{O}(\eps^{N+3})$--bound is in the $C^\infty$-sense. This implies that
\begin{multline}
p_\eps \circ \exp(\eps H_{G_{\leq N+1}}) \\
= p + \eps \left\langle q \right\rangle + \eps^2 q_2 + \ldots + \eps^{N+1}q_{N+1} + (r_{N+2} - H_p G_{N+1})\eps^{N+2} +
\mathcal{O}(\eps^{N+3}).
\end{multline}
As above, we may find a smooth real-valued solution of $H_p G_{N+1} = r_{N+2} - \left\langle r_{N+2} \right\rangle$, defined near $p^{-1}(0)$.

\medskip
\noindent
The functions $G_j$, $j \geq 0$, are defined in a fixed neighborhood of $p^{-1}(0)$, and by Borel's lemma we may choose $G(x,\xi,\eps) \in C^{\infty}$
near $p^{-1}(0)$, smooth in $\eps \in {\rm neigh}(0,\real)$, which is given by
\begeq
\label{eq2.7}
G \sim \sum_{j=0}^\infty \eps^j G_j,
\endeq
asymptotically in the $C^{\infty}$--sense. We have then achieved that $p_\eps \circ \exp(\eps H_G)$ is in involution with $p$ modulo
$\mathcal{O}(\eps^\infty)$, in a fixed neighborhood of $p^{-1}(0)$, as desired.

\bigskip
\noindent
Now an application of Cartan's formula shows that the canonical transformation $\exp(\eps H_G)$ is exact in the sense that the 1-form
$(\exp(\eps H_G))^*\lambda - \lambda$ is exact, where $\lambda$ is the fundamental 1-form on $T^*M$. By Egorov's theorem, we may therefore quantize the
real-valued smooth canonical transformation $\exp(\eps H_G)$ by a (single-valued) $h$-Fourier integral operator
$U_{\eps}= {\cal O}(1): L^2(M)\to L^2(M)$, which is microlocally unitary near $p^{-1}(0)$. Then we have that the selfadjoint operator
$\widetilde P_\eps := U_{\eps}^{-1}P_\eps U_{\eps}$, defined microlocally near $p^{-1}(0)$, is such that its leading symbol is of the form
$p + \eps \left\langle q \right\rangle + \mathcal{O}(\eps^2)$, where the $\mathcal{O}(\eps^2)$ term is in involution with
$p$, modulo $\mathcal{O}(\eps^\infty)$.

\medskip
\noindent
Furthermore, by the results of Section 2 of \cite{HS1}, we know that if we choose the principal symbol of the microlocally unitary Fourier integral operator
${U}_{\eps}$ to be of constant argument, then ${U}_{\eps}$ enjoys the improved Egorov property, so that on the level of symbols we have
$\widetilde P_\eps = P_\eps \circ \exp(\eps H_G) + \mathcal{O}(h^2)$. A natural choice of $U_{\eps}$ is therefore given by
$U_{\eps} = e^{-i\eps G/h}$, since then the principal symbol of $U_{\eps}$ solves a real transport equation, using also that the subprincipal
symbol of $G$ vanishes.

\bigskip
\noindent
We summarize the discussion above in the following result.
\medskip
\noindent
\begin{prop}
There exists $G(x,\xi,\eps)\in C^{\infty}({\rm neigh}(p^{-1}(0), T^*M))$ real-valued, depending smoothly on $\eps\in {\rm neigh}(0,\real)$, with the asymptotic
expansion {\rm (\ref{eq2.7})} in the space of real-valued smooth functions in a fixed \neigh{} of $p^{-1}(0)$, such that the microlocally defined selfadjoint
operator
$$
\widetilde P_\eps = e^{i\eps G/h} P_{\eps} e^{-i\eps G/h}
$$
depends on $\eps$ in a $C^{\infty}$-fashion and has the leading symbol of the form
$p + \eps \left\langle q \right\rangle + \mathcal{O}(\eps^2)$, where the ${\cal O}(\eps^2)$--term Poisson commutes with $p$ modulo
${\cal O}(\eps^{\infty})$. The subprincipal symbol of $\widetilde P_\eps$ is ${\cal O}(\eps)$. Assume furthermore that $\eps \ll h$, so that the spectrum of $P_{\eps}$ near $0$ retains a cluster structure, being confined to the union of intervals
$$
I_k(\eps) = f(h(k-\theta)) + [-{\cal O}(h^2 + \eps), {\cal O}(h^2 + \eps)],\quad k\in \z.
$$
If $z\in {\rm Spec}(P_{\eps})\cap {\rm neigh}(0,\real)$ is such that $z\in I_k(\eps)$, for some $k$, then we have
$$
\eps {\rm min}_{{\rm neigh}\, (p^{-1}(0))} \langle{q\rangle} -{\cal O}(\eps^2 + h^2) \leq z - f(h(k-\theta)) \leq
\eps {\rm max}_{{\rm neigh}\,(p^{-1}(0))} \langle{q\rangle} + {\cal O}(\eps^2 + h^2).
$$
\end{prop}

\noindent
Here the last estimate follows by an application of sharp G\aa{}rding's inequality.

\section{Normal form near a Lagrangian torus}
\setcounter{equation}{0}
In Proposition 2.2, we have reduced ourselves to a microlocally defined selfadjoint operator $\widetilde{P}_{\epsilon}$, acting on $L^2(M)$, with the leading
symbol of the form
$$
p + \eps \langle{q\rangle} + {\cal O}(\eps^2),
$$
where the ${\cal O}(\eps^2)$--term Poisson commutes with $p$, modulo ${\cal O}(\eps^{\infty})$. The subprincipal symbol of $\widetilde{P}_{\eps}$ is
${\cal O}(\eps)$. In what follows, when working with the operator $\widetilde{P}_{\eps}$, to simplify the notation, we shall drop the tilde
and write $P_{\eps}$ instead.

\medskip
\noindent
Let $F_0 \in \real$ be such that ${\rm min}_{p^{-1}(0)} \langle{q\rangle} < F_0 < {\rm max}_{p^{-1}(0)} \langle{q\rangle}$ and assume that $F_0$ is a regular
value of $\langle{q\rangle}$, viewed as a function on the space of closed orbits $\Sigma$. After replacing $q$ by $q-F_0$ we may assume that $F_0 = 0$, and let us
consider the $H_p$--flow invariant set
$$
\Lambda_{0,0} = \{\rho\in T^*M;\,\, p(\rho) = 0, \left\langle q \right\rangle(\rho)= 0\}.
$$
We know that $T(0)$ is the minimal period of all closed $H_p$--trajectories in $\Lambda_{0,0}$ and since $dp, d\left\langle q \right\rangle$ are linearly
independent at each point of $\Lambda_{0,0}$, we see that $\Lambda_{0,0}$ is a Lagrangian manifold which is a union of finitely many tori. Assume for simplicity
that $\Lambda_{0,0}$ is connected so that it is equal to a single Lagrangian torus. Since the functions $p,\left\langle q \right\rangle$ are in involution, they
form a completely integrable system in a neighborhood of $\Lambda_{0,0}$. We have action-angle coordinates near $\Lambda_{0,0}$~\cite{HZ}, given by a smooth
real-valued canonical transformation
\begeq
\label{eq3.1}
\kappa: {\rm neigh}(\xi = 0,T^*\mathbf{T}^2)\to {\rm neigh}(\Lambda_{0,0},T^*M),\quad {\bf T}^2 = \real^2/2\pi \z^2,
\endeq
mapping the zero section in $T^*{\bf T}^2$ onto $\Lambda_{0,0}$, and such that $p\circ \kappa = p(\xi)$,
$\left\langle q \right\rangle\circ \kappa = \left\langle q \right\rangle(\xi)$. Here we make the identification
$T^*\mathbf{T}^2 \cong (\mathbf{R}/2\pi\mathbf{Z})^2_x \times \mathbf{R}^2_\xi$. Since the classical flow of $p$ is periodic with minimal period $T(0)$ in
$\Lambda_{0,0}$, we may and will choose $\kappa$ so that in fact $p\circ \kappa = p(\xi_1)$, by letting $\xi_1$ be the normalized action of a closed
$H_p$-trajectory of minimal period --- see the discussion in Section 4 of \cite{HS1}. The linear independence of the differentials of $p$ and
$\langle{q\rangle}$ implies that $p'(0) \neq 0$, $\partial_{\xi_2}\left\langle q \right\rangle(0) \neq 0$. We may also remark that when expressed in terms of the
action coordinate $\xi_1$, the function $p$ becomes $p(\xi_1) = f(\xi_1)$, where the smooth function $f$ has been introduced after (\ref{eq1.85}).

\bigskip
\noindent
Implementing the real canonical transformation $\kappa$ in (\ref{eq3.1}) by means of a multi-valued microlocally unitary $h$-Fourier integral operator
$U: L^2_f({\bf T}^2) \rightarrow L^2(M)$, which also has the improved Egorov property~\cite{HS1}, we get a new selfadjoint operator $U^{-1} P_{\eps} U$,
which for simplicity, will still be denoted by $P_{\eps}$,
\begeq
\label{eq3.2}
P_\eps: L^2_f(\mathbf{T}^2)\to L^2_f(\mathbf{T}^2).
\endeq
Here the operator $P_{\eps}$ is defined microlocally near $\xi=0$ in $T^*{\bf T}^2$, with the full (Weyl) symbol of the form
\begeq
\label{eq3.3}
P_{\eps} \sim \sum_{j=0}^\infty h^j p_j(x,\xi,\eps),
\endeq
the principal symbol being
\begeq
\label{eq3.4}
p_0(x,\xi,\eps) = p(\xi_1) + \eps\left\langle q \right\rangle(\xi) + \mathcal{O}(\eps^2)
\endeq
with the $\mathcal{O}(\eps^2)$ error term independent of $x_1$ modulo $\mathcal{O}(\eps^\infty)$. Furthermore,
$$
p_1(x,\xi,\eps) = \mathcal{O}(\eps),
$$
and all the terms in the expansion (\ref{eq3.3}) are smooth and  real-valued. The dependence on $\eps\in [0,\eps_0)$
in (\ref{eq3.3}) is still $C^{\infty}$. The space $L^2_f(\mathbf{T}^2)$ here stands for the subspace of $L^2_{\rm loc}(\real^2)$ consisting of Floquet periodic
functions $u(x)$, satisfying
$$
u(x - \nu) = e^{i\nu\cdot \Theta} u(x),\quad \nu \in (2\pi \mathbf{Z})^2, \quad \Theta = \frac{S}{2\pi h} + \frac{\alpha}{4}.
$$
Here $S=(S_1,S_2)$ with $S_j$ being the action of the generator $\gamma_j$ of the first homotopy group of $\Lambda_{0,0}$, with $\gamma_1$ being given by a closed
$H_p$--trajectory of minimal period, and $\alpha=(\alpha_1,\alpha_2)$ is the corresponding Maslov index.

\subsection{Removing the $x_1$ dependence}
Our first goal is to eliminate the $x_1$-dependence in $p_j$, $j \geq 1$, in (\ref{eq3.3}). Let $A=A(x,\xi,\eps)\in C^{\infty}$ be real-valued, $x\in {\bf T}^2$,
$\xi \in {\rm neigh}(0,\real^2)$, smooth in $\eps \geq 0$, and let us consider the conjugation of the selfadjoint operator $P_{\eps}$ by the unitary
$h$-pseudodifferential operator $e^{iA^w}$. We have
$$
e^{-iA^w}P_\eps e^{iA^w} = e^{-i{\rm ad} A^w} P_{\eps} = \sum_{k=0}^{\infty} \frac{(-i{\rm ad} A^w)^k}{k!} P_{\eps},
\quad ({\rm ad} A^w) P_{\eps} = [A^w,P_{\eps}].
$$
Identifying the symbols with the corresponding $h$-Weyl quantizations, we obtain that
\begin{align*}
e^{-iA^w}P_\eps e^{iA^w} &= P_\eps  + e^{-iA^w}[P_\eps,e^{iA^w}] \\
&= p_0(x,\xi,\eps) + h(p_1(x,\xi,\eps) + H_{p_0}A(x,\xi,\eps)) + \mathcal{O}(h^2).
\end{align*}
We shall now show that $A$ can be chosen real-valued smooth, so that $p_1 + H_{p_0} A$ becomes independent of $x_1$, modulo ${\cal O}(\eps^{\infty})$. In doing so,
we shall construct the $C^{\infty}$--symbol $A$ as a formal power series in $\eps$. Introducing the Taylor expansions,
$$
p_0(x,\xi,\eps) \sim \sum_{\ell =0}^{\infty} \eps^{\ell} p_{0,\ell}(x,\xi),\quad p_1(x,\xi,\eps) \sim \sum_{\ell=1}^{\infty} \eps^{\ell} p_{1,\ell}(x,\xi),
$$
and writing
$$
A\sim \sum_{\ell =1}^{\infty} \eps^{\ell} a_{\ell}(x,\xi),
$$
we compute the power series expansion of the Poisson bracket,
$$
H_{p_0}A \sim \sum_{k\geq 0, \ell \geq 1} \eps^{k+\ell} \{p_{0,k},a_{\ell}\} = \sum_{m=1}^{\infty} \eps^m f_m,
$$
where
$$
f_m = \sum_{k+\ell=m, \atop k\geq 0, \ell \geq 1} \{p_{0,k},a_{\ell}\}.
$$
We would like to choose the coefficients $a_{\ell}$, $\ell \geq 1$, so that $p_{1,\ell} + f_{\ell}$ is independent of $x_1$, for all $\ell \geq 1$. When $\ell=1$,
we have $p_{1,1} + f_1 = p_{1,1} + \partial_{\xi_1}p\, \partial_{x_1} a_1$, and since $\partial_{\xi_1}p(0)\neq 0$, we can determine $a_1$ real-valued
by solving the transport equation,
$$
p_{1,1} + \partial_{\xi_1}p\, \partial_{x_1} a_1 = \langle{p_{1,1}\rangle}_{x_1},
$$
the right hand side standing for the average of $p_{1,1}$ with respect to $x_1$. Arguing inductively, assume that the smooth real-valued functions
$a_1,\ldots\, a_m$ have already been determined. The term $p_{1,m+1} + f_{m+1}$ is of the form
$$
p_{1,m+1} + \partial_{\xi_1} p\, \partial_{x_1} a_{m+1} + \sum_{k+\ell = m+1,\, \ell < m+1} \{p_{0,k},a_{\ell}\},
$$
and it is therefore clear that we can choose $a_{m+1}$ real, so that this expression becomes independent of $x_1$. Arguing in this fashion, we obtain a sequence
$a_j \in C^{\infty}({\rm neigh}(\xi=0, T^*{\bf T}^2))$, $a_j$ real-valued, so that if $A\in C^{\infty}$ in all variables and real-valued, is such that
$$
A(x,\xi,\eps) \sim \sum_{j=1}^{\infty} \eps^j a_j,
$$
then the subprincipal symbol of the selfadjoint operator $e^{-iA^w} P_{\eps} e^{iA^w}$ is ${\cal O}(\eps)$ and independent of $x_1$, modulo
${\cal O}(\eps^{\infty})$.

\medskip
\noindent
Assume inductively that we have found $A_0=A,\ldots\, A_{N-1}$ real-valued so that the selfadjoint operator
$$
P_{\eps}^{(N)}:= e^{-i{\rm ad}(h^{N-1} A_{N-1})} \circ \ldots e^{-i{\rm ad}(h A_1)}\circ e^{-i{\rm ad}(A)} P_{\eps}
$$
is of the form $\sim \sum_{j=0}^{\infty} h^j p_j$, where $p_j$ are independent of $x_1$ modulo ${\cal O}(\eps^{\infty})$, for $j\leq N$. We then look for
a conjugation by a unitary $h$-pseudodifferential operator of the form $e^{ih^N A_N^w}$, and we see as before that the leading symbol of
$e^{-ih^N A_N^w} [P_{\varepsilon}^{(N)},e^{ih^N A_N^w}]$ is $h^{N+1} H_{p_0} A_N$.
Therefore, $e^{-ih^N A_N^w} P_{\varepsilon}^{(N)} e^{ih^N A_N^w}$ is of the form $\sim \sum_{j=0}^{\infty} h^j \widetilde{p}_j$, where $\widetilde{p}_j = p_j$
for $j\leq N$, and $\widetilde{p}_{N+1} = p_{N+1} + H_{p_0} A_N$. It is therefore clear that we can determine $A_N$, real-valued and smooth, as a formal power
series in $\eps$, so that $\widetilde{p}_{N+1}$ becomes independent of $x_1$, modulo ${\cal O}(\eps^{\infty})$. Using the Campbell-Hausdorff formula
and Borel's lemma, we see that there exists a selfadjoint $h$-pseudodifferential operator ${\cal A}$ with a real-valued symbol
$\sim \sum_{\nu=0}^{\infty} h^{\nu} a_{\nu}(x,\xi,\eps) \in S(1)$, smooth in $\eps \in [0,\eps_0)$, with $a_0 = {\cal O}(\eps)$, , such that
$$
e^{-i{\rm ad} {\cal A}^w} \sim \ldots e^{-ih^2 {\rm ad}(A_2^w)} \circ e^{-ih {\rm ad}(A_1^w)} \circ e^{-i {\rm ad}(A_0^w)},
$$
and we conclude that the selfadjoint operator $\widetilde{P}_{\eps} = e^{-i{\rm ad} {\cal A}^w} P_{\eps}$ is of the form
$$
\widetilde{P}_{\eps} = \sum_{j=0}^{\infty} h^j \widetilde{p}_j(x_2,\xi,\eps),
$$
modulo ${\cal O}((h,\eps)^{\infty})$. Here $\widetilde{p}_j$ are real-valued, smooth in $\eps \in [0,\eps_0)$ and independent of $x_1$, with
$\widetilde{p}_0(x_2,\xi,\eps) = p(\xi_1) + \eps \langle{q\rangle}(\xi) + {\cal O}(\eps^2)$, and $\widetilde{p}_1(x_2,\xi,\eps) = {\cal O}(\eps)$.

\subsection{Removing the $x_2$ dependence}
In the previous subsection, using only the fact that $\partial_{\xi_1}p(0)\neq 0$, we have carried out repeated averagings along the $H_p$--flow, thus
eliminating the $x_1$-dependence in the full symbol of $P_{\eps}$ by means of a conjugation by an elliptic unitary $h$-pseudodifferential operator. We have
reduced ourselves to a selfadjoint operator of the form,
\begeq
\label{eq3.5}
\widetilde P_\eps \sim \sum_{j=0}^\infty h^j \widetilde p_j(x_2,\xi,\eps)\quad \text{on }L^2_f(\mathbf{T}^2),
\endeq
modulo ${\cal O}((\eps,h)^{\infty})$, where $\widetilde{p}_j$ in (\ref{eq3.5}) are real-valued, smooth in $\eps \in [0,\eps_0)$, with
$$
\widetilde p_0(x_2,\xi,\eps) = p(\xi_1) + \eps \left\langle q \right\rangle(\xi) + \mathcal{O}(\eps^2)
$$
with the $\mathcal{O}(\eps^2)$ error term independent of $x_1$. Also, $\widetilde{p}_1(x_2,\xi,0) = 0$.

\bigskip
\noindent
Continuing to argue in the spirit of \cite{HS1}, \cite{HS2}, we shall now look for an additional conjugation by means of unitary $h$--Fourier integral operators,
which eliminates the $x_2$-dependence in the full symbol in (\ref{eq3.5}). Following \cite{HS1}, to that end it will be convenient to construct the conjugating operator by
viewing $\eps$ and $h^2/\eps$ as two independent small parameters, provided that $\eps$ is not too small.

\medskip
\noindent
On the level of symbols, we write, using that $\widetilde{p}_1(x_2,\xi,\eps) = \eps q_1(x_2,\xi,\eps)$, where $q_1$ is real-valued and $C^{\infty}$
in all variables,
\begin{multline}
\label{eq3.51}
\widetilde{P}_\eps = p(\xi _1)+\eps \left(\langle q\rangle
(\xi )+{\cal O}(\eps) + h q_1(x_2,\xi ,\eps)+\frac{h^2}{\eps} \widetilde{p}_2 + h\frac{h^2}{\eps}\widetilde{p}_3+\ldots\right)\\
 = p(\xi _1)+\eps \left (r_0(x_2,\xi,\eps,\frac{h^2}{\eps}) + hr_1(x_2,\xi ,\eps,\frac{h^2}{\eps}) + h^2 r_2+\ldots\right),
\end{multline}
with
\begeq
\label{eq3.52}
r_0(x_2,\xi,\eps,\frac{h^2}{\eps}) = \langle{q\rangle}(\xi) + {\cal O}(\eps) + \frac{h^2}{\eps} \widetilde{p}_2,
\endeq
$$
r_1(x_2,\xi,\eps,\frac{h^2}{\eps}) = q_1(x_2,\xi,\eps) + \frac{h^2}{\eps} \widetilde{p}_3,
$$
$$
r_j(x_2,\xi,\eps,\frac{h^2}{\eps}) = \frac{h^2}{\eps} \widetilde{p}_{j+2},\quad j\geq 2.
$$

\medskip
\noindent
When eliminating the variable $x_2$, let us introduce the basic assumption that
\begeq
\label{eq3.6}
\eps = {\cal O}(h^{\delta}),\quad \delta>0,
\endeq
and also, assume that
\begeq
\label{eq3.7}
\frac{h^2}{\eps}\leq \delta_0,
\endeq
for some $\delta_0 >0$ sufficiently small but independent of $h$. Replacing first (\ref{eq3.7}) by the strengthened hypothesis,
\begeq
\label{eq3.8}
\frac{h^2}{\eps} \leq {\cal O}(h^{\delta_1}),\quad \delta_1 > 0,
\endeq
let us describe the construction of a unitary conjugation elimi\-na\-ting the $x_2$--depen\-dence in $\widetilde{P}_{\eps}$.

\medskip
\noindent
When $b_0 = b_0(x_2,\xi,\eps, \frac{h^2}{\eps})$ is real-valued and smooth for $\xi \in {\rm neigh}(0,\real^2)$, $\eps$,
$h^2/\eps\in [0,\eps_0)$, and is such that $b_0 = {\cal O}(\eps + h^2/\eps)$ in the $C^{\infty}$--sense, we consider the selfadjoint
operator
\begeq
\label{eq3.9}
e^{\frac{i}{h} B_0} \widetilde{P}_{\epsilon} e^{-\frac{i}{h}B_0},\quad B_0 = b_0^w(x_2,hD_x,\eps, h^2/\eps).
\endeq
Since the commutator $[B_0,p(hD_{x_1})] = 0$, we see that the full symbol of the conjugated operator (\ref{eq3.9}) is real-valued and of the form
$$
p(\xi_1) + \eps \left(\widehat{r}_0 + h\widehat{r}_1 + \ldots\right),
$$
where by Egorov's theorem,
$$
\widehat{r}_0 = r_0 \circ \exp(H_{b_0}) = \sum_{k=0}^{\infty} \frac{1}{k!} H_{b_0}^k r_0,
$$
while $\widehat{r}_j = {\cal O}(1)$ for $j\geq 1$. Since the canonical transformation $\exp(H_{b_0})$ is exact, we see that the conjugated operator still acts
on the space $L^2_f({\bf T}^2)$ of Floquet periodic functions.

\medskip
\noindent
It follows from (\ref{eq3.52}) that
$$
\widehat{r}_0 = \langle{q\rangle}(\xi) + {\cal O}\left(\eps + \frac{h^2}{\eps}\right) - \partial_{\xi_2}\langle{q\rangle} \partial_{x_2} b_0 +
{\cal O}\left(\left(\eps,\frac{h^2}{\eps}\right)^2\right),
$$
and using that $\partial_{\xi_2}\langle{q\rangle} \neq 0$, it becomes clear how to construct a real-valued smooth symbol
$b_0 = {\cal O}(\eps + h^2/\eps)$, defined near $\xi=0$ in $T^*{\bf T}^2$, smooth in $\eps$, $h^2/\eps \in {\rm neigh}(0,\real)$, as a
formal Taylor series in $\eps$, $h^2/\eps$, so that $\widehat{r}_0 = \langle{q\rangle} + {\cal O}(\eps + h^2/\eps)$ is independent of $x$, modulo
${\cal O}(h^{\infty})$, in view of (\ref{eq3.6}), (\ref{eq3.8}).

\bigskip
\noindent
Dropping the assumption (\ref{eq3.8}), we now come to discuss the construction of a conjugating Fourier integral operator when only (\ref{eq3.7}) is valid.
Following~\cite{HS1}, we consider the eikonal equation for $\varphi = \varphi(x_2,\xi,\eps,h^2/\eps)$,
\begeq
\label{eq3.10}
r_0 \left(x_2,\xi_1, \xi_2 + \partial_{x_2} \varphi, \eps, \frac{h^2}{\eps}\right) = \langle{r_0(\cdot,\xi, \eps,\frac{h^2}{\eps})\rangle},
\endeq
where $\langle{\cdot\rangle}$ in the right hand side stands for the average with respect to $x_2$. Since $\partial_{\xi_2}\langle{q\rangle} \neq 0$, by
Hamilton-Jacobi theory, (\ref{eq3.10}) has a smooth real-valued solution with $\partial_{x_2}\varphi$ single-valued and
$$
\partial_{x_2} \varphi = {\cal O}\left(\eps + \frac{h^2}{\eps}\right).
$$
Taylor expanding (\ref{eq3.10}) and using that
$$
\partial_{\xi_2} r_0 \left(x_2,\xi,\eps,\frac{h^2}{\eps}\right) = \partial_{\xi_2}\langle{q\rangle}(\xi) + {\cal O}\left(\eps + \frac{h^2}{\eps}\right),
$$
we get
$$
\varphi = \varphi_{\rm per} + x_2\zeta_2,
$$
where $\varphi_{\rm per} = {\cal O}(\eps + h^2/\eps)$ is periodic in $x_2$ and $\zeta_2 = \zeta_2(\xi,\eps,h^2/\eps) = {\cal O}((\eps, h^2/\eps)^2)$. Let us set
$$
\eta = \eta\left(\xi,\eps,\frac{h^2}{\eps}\right) = (\xi_1, \xi_2 + \zeta_2),
$$
and
$$
\psi\left(x_2,\eta, \eps,\frac{h^2}{\eps}\right) = \varphi_{{\rm per}} + x\cdot \eta,
$$
where $\varphi_{\rm per}$ is viewed as a function of $\eta$ rather than $\xi$. Associated to the function $\psi$ is the real-valued smooth canonical transformation
\begeq
\label{eq3.11}
\kappa: (\psi'_{\eta}, \eta) \to (x,\psi'_x),
\endeq
which is an ${\cal O}(\eps + h^2/\eps)$--perturbation of the identity in the $C^{\infty}$--sense, and such that if
$(x,\xi) = \kappa(y,\eta)$, then $\xi_1 = \eta_1$. We have by construction,
\begin{multline*}
\left(r_0 \circ \kappa\right)(y,\eta,\eps, h^2\eps) = r_0(x,\psi'_x,\eps, h^2/\eps) = \langle{r_0(\cdot, \xi, \eps,\frac{h^2}{\eps})\rangle} \\
=  \langle{r_0(\cdot, \eta, \eps,\frac{h^2}{\eps})\rangle} + {\cal O}\left(\left(\eps, h^2/\eps\right)^2\right),
\end{multline*}
which is a function of $(y,\eta)$, independent of $y$.

\medskip
\noindent
We can quantize the canonical transformation $\kappa$ in (\ref{eq3.11}) by a microlocally unitary Fourier integral operator, and after conjugation by this operator, we obtain a new operator, still denoted by
$\widetilde{P}_{\eps}$, which is of the form (\ref{eq3.51}), where
$$
r_0 = \langle{q\rangle}(\xi) + {\cal O}\left(\eps + \frac{h^2}{\eps}\right)
$$
is independent of $x$, and $r_j = {\cal O}(1)$ in the $C^{\infty}$--sense, for $j\geq 1$. Furthermore, as explained in Section 4 of~\cite{HS1}, the conjugated
operator $\widetilde{P}_{\eps}$ still acts on the space $L^2_f({\bf T}^2)$ of Floquet periodic functions.

\bigskip
\noindent
Let us consider therefore an operator of the form
$$
\widetilde{P}_{\eps} = p(\xi_1) + \eps\left(r_0(\xi,\eps, \frac{h^2}{\eps}) + hr_1(x_2,\xi,\eps,\frac{h^2}{\eps}) +\ldots\right),
$$
where $r_0 = \langle{q\rangle}(\xi) + {\cal O}(\eps + h^2/\eps)$ is independent of $x$, and $r_j = {\cal O}(1)$, $j\geq 1$. Furthermore, all the
terms $r_j$ are real-valued, smooth, and depend smoothly on $\eps$, $h^2/\eps \in {\rm neigh}(0,\real)$. To eliminate the $x_2$--dependence in the lower order
terms $r_j$, $j\geq 1$, we could argue as in the previous step, making the terms $r_j$ independent of $x_2$ one at a time, but here we would like to
describe a slightly different method, which has the merit of being more direct. Let us look for a conjugation by an elliptic unitary pseudodifferential operator
of the form $e^{iB/h}$, where
$$
B(x_2,\xi,\eps,\frac{h^2}{\eps};h) = \sum_{\nu=1}^{\infty} h^{\nu} b_{\nu} (x_2,\xi,\eps, \frac{h^2}{\eps}).
$$
Here $b_{\nu}$ are real-valued smooth and depend smoothly on $\eps$, $h^2/\eps\in {\rm neigh}(0,\real)$.
The conjugated operator
$$
e^{\frac{i}{h} B} \widetilde{P}_{\eps} e^{-\frac{i}{h}B} = e^{\frac{i}{h}{\rm ad}\, B} \widetilde{P}_{\epsilon} = \sum_{k=0}^{\infty}
\frac{(i {\rm ad}\, B)^k}{h^k k!} \widetilde{P}_{\epsilon}
$$
is selfadjoint and can be expanded as follows,
\begin{multline}
\label{symbol4}
p(\xi_1) + \eps \sum_{k=0}^\infty \sum_{j_1=1}^\infty...\sum_{j_k=1}^\infty \sum_{\ell =0}^\infty
h^{\ell + j_1+..+j_k} \frac{1}{k!} \left(\frac{i}{h}{\rm ad}\, b_{j_1} \right)..\left (\frac{i}{h}{\rm ad}\, b_{j_k}\right) r_\ell \\
= p(\xi _1)+\eps \sum_{n=0}^ \infty h^n\widehat{r}_n.
\end{multline}
Here $\widehat{r}_n$ is equal to the sum of all the coefficients for $h^n$ coming from the expressions
\begin{equation}
\label{symbol5}
h^{\ell + j_1+..+j_k} \frac{1}{k!} \left(\frac{i}{h}{\rm ad}\, b_{j_1} \right)..\left (\frac{i}{h}{\rm ad}\, b_{j_k}\right) r_\ell,
\end{equation}
with $\ell+j_1+..+j_k\le n$ and $j_\nu \ge 1$. In particular, we see that $\widehat{r}_n$ are all real-valued, thanks to the observation that if $A$, $B$ are
selfadjoint, then so is the operator $i[A,B] = (i{\rm ad}A)B$. Then $\widehat{r}_0=r_0$,
$\widehat{r}_1=r_1+H_{b_1}r_0=r_1-H_{r_0}b_1$,..,$\widehat{r}_n=r_n-H_{r_0}b_n +s_n$, where $s_n$ only depends on
$b_1,...,b_{n-1}$ and is the sum of all coefficients of $h^n$ arising in the expressions (\ref{symbol5}) with $\ell +j_1 +..+j_k\le n$,
$j_1,..,j_k,\ell <n$, $j_\nu \ge 1$.

\medskip
\noindent
It is therefore clear how to find $b_1,b_2,\ldots $ real-valued smooth, successively, with $b_j={\cal O}(1)$, such that all the coefficients $\widehat{r}_j$ in
(\ref{symbol4}) are independent of $x$ and $={\cal O}(1)$.

\bigskip
\noindent
The discussion in this section may be summarized in the following theorem.

\begin{theo}
Let us make all the general assumptions of Section {\rm 1} and let $F_0 \in \real$ be a regular value of $\langle{q\rangle}$, viewed as a function on the
space of closed orbits $\Sigma$. Assume that the Lagrangian manifold
$$
\Lambda_{0,F_0}: p=0,\quad \langle{q\rangle} = F_0
$$
is connected. When $\gamma_1$ and $\gamma_2$ are the fundamental cycles in $\Lambda_{0,F_0}$ with $\gamma_1$ corresponding to a closed $H_p$--trajectory of
minimal period, we write $S = (S_1,S_2)$ and $\alpha= (\alpha_1,\alpha_2)$ for the actions and the Maslov indices of the cycles, respectively. Assume furthermore
that $\eps = {\cal O}(h^{\delta})$, $\delta >0$, is such that $h^2/\eps \leq \delta_0$, for some $\delta_0 >0$ sufficiently small but fixed. There exists a
smooth Lagrangian torus $\widehat{\Lambda}_{0,F_0}\subset T^*M$, which is an ${\cal O}(\eps + h^2/\eps)$--perturbation of $\Lambda_{0,F_0}$ in the
$C^{\infty}$--sense, such that when $\rho\in T^*M$ is away from a small neighborhood of $\widehat{\Lambda}_{0,F_0}$ and $\abs{p(\rho)} \leq 1/C$, for $C>0$
sufficiently large, we have
$$
\abs{\langle{q\rangle}(\rho) - F_0} \geq \frac{1}{{\cal O}(1)}.
$$
There exists a $C^{\infty}$ real-valued canonical transformation
$$
\kappa: {\rm neigh}(\widehat{\Lambda}_{0,F_0},T^*M) \to {\rm neigh}(\xi =0,T^*{\bf T}^2),
$$
mapping to $\widehat{\Lambda}_{0,F_0}$ to $\xi =0$, and a corresponding uniformly bounded $h$-Fourier integral operator
$$
U = {\cal O}(1): L^2(M) \rightarrow L^2_f({\bf T}^2),
$$
which has the following properties:
\begin{enumerate}
\item The operator $U$ is microlocally unitary near $\widehat{\Lambda}_{0,F_0}$: if
$U^* = {\cal O}(1): L^2_{f}({\bf T}^2) \rightarrow L^2(M)$ is the complex adjoint, then for every
$\chi_1 \in C_0^{\infty}({\rm neigh}(\widehat{\Lambda}_{0,F_0},T^*M))$, we have
\begin{equation}
\left(U^* U - 1\right){\rm Op}_h^w(\chi_1) = {\cal O}(h^{\infty}): L^2(M) \rightarrow L^2(M).
\end{equation}
For every $\chi_2 \in C_0^{\infty}({\rm neigh}(\xi=0, T^*{\bf T}^2))$, we have
$$
\left(U U^*- 1\right)\chi_2^w(x,hD_x)={\cal O}(h^{\infty}): L^2_{f}({\bf T}^2)\rightarrow L^2_{f}({\bf T}^2).
$$
\item We have a normal form for $P_{\eps}$: Acting on $L^2_{f}({\bf T}^2)$, there exists a selfadjoint operator
$\widehat{P}\left(hD_x,\eps,\frac{h^2}{\eps};h\right)$ with the symbol
$$
\widehat{P}\left(\xi,\eps,\frac{h^2}{\eps};h\right) \sim p(\xi_1)+\eps \sum_{j=0}^{\infty} h^j r_j\left(\xi,\eps,\frac{h^2}{\eps}\right),
\quad \abs{\xi}\leq \frac{1}{{\cal O}(1)},
$$
smooth in $\xi \in {\rm neigh}(0,\real^2)$, and smooth in $\eps$, $h^2/\eps \in {\rm neigh}(0,\real)$, such that
$$
r_0=\langle{q}\rangle(\xi)+{\cal O}\left(\eps +\frac{h^2}{\eps}\right),
$$
and
$$
r_j={\cal O}(1),\quad j\geq 1,
$$
and such that $\widehat{P}U = U P_{\eps}$ microlocally near $\widehat{\Lambda}_{0,F_0}$, i.e.
$$
\left(\widehat{P}U - U P_{\eps}\right){\rm Op}^w_h(\chi_1) = {\cal O}(h^{\infty}),
\quad \chi_2^w(x,hD_x)\left(\widehat{P}U - U P_{\epsilon}\right) = {\cal O}(h^{\infty}),
$$
in the operator sense, for every $\chi_1$, $\chi_2$ as in {\rm 1)}.
\end{enumerate}
\end{theo}

\section{Eigenvalue asymptotics in subclusters}
\setcounter{equation}{0}
Throughout this section, we shall assume that $\eps \ll h$ and that the lower bound $h^2/\eps \leq \delta_0 \ll 1$ is valid. We then know that Theorem 3.1 applies
and that the spectrum of $P_{\eps}$ near $0$ is confined to the union of intervals,
$$
I_k(\eps) = f(h(k-\theta)) + [-{\cal O}(\eps), {\cal O}(\eps)],\quad k\in \z,\quad \theta = \frac{S_1}{2\pi h} + \frac{\alpha_1}{4},
$$
disjoint for all $h>0$ small enough.

\bigskip
\noindent
When proving Theorem 1.2, following~\cite{HS2}, let us first check that if $z\in {\rm neigh}(0,\real)$ is such that
\begeq
\label{eq4.1}
\abs{z-f(h(k -\theta)) - \eps F_0} \leq \frac{\eps}{C},\quad C\gg 1,
\endeq
for some $k \in \z$, and $z$ avoids the union of the pairwise disjoint open intervals $J_{\ell}(h)$ of length $\eps h/{\cal O}(1)$,
that are centered at the quasi--eigenvalues
\begeq
\label{eq4.1.1}
\widehat{P}
\left(h(k -\frac{\alpha_1}{4})-\frac{S_1}{2\pi},h(\ell -\frac{\alpha_2}{4})-\frac{S_2}{2\pi},\eps,\frac{h^2}{\eps};h\right),
\endeq
for $\ell \in \z$, then the operator
$$
P_{\eps}-z: H(m) \rightarrow L^2(M)
$$
is bijective.

\medskip
\noindent
To that end, consider a partition of unity on $T^*M$,
\begeq
\label{eq4.2}
1=\chi+\psi_{1,+}+\psi_{1,-}+\psi_{2,+}+\psi_{2,-}.
\endeq
Here $\chi\in C^{\infty}_0(T^*M)$ is supported in a small flow invariant neighborhood of $\widehat{\Lambda}_{0,F_0}$ where the operator $U$ of Theorem 3.1
is defined and unitary, and where $P_{\eps}$ is intertwined with $\widehat{P}$, and $\chi=1$ near $\widehat{\Lambda}_{0,F_0}$. Thanks to Theorem 3.1, we also assume, as we may, that on
the operator level,
\begeq
\label{eq4.3}
[P_{\eps}, \chi]={\cal O}(h^{\infty}): L^2 \rightarrow L^2.
\endeq
Furthermore, the functions $\psi_{1,\pm}\in C^{\infty}_0(T^*M)$ are supported in flow invariant regions $\Omega_{\pm}$, such that
$\pm (\langle{q}\rangle - F_0) \geq 1/{\cal O}(1)$ in $\Omega_{\pm}$, respectively. Moreover, we can arrange so that $\psi_{1,\pm}$ are in involution with
$p$, the principal symbol of $P_{\eps=0}$. Finally, $\psi_{2,\pm}\in C^{\infty}_b(T^*M)$ are such that $\pm p >1/{\cal O}(1)$ in the support of
$\psi_{2,\pm}$.

\medskip
\noindent
Let us consider the equation,
$$
(P_{\eps}-z)u = v,\quad u\in H(m),
$$
when $z\in {\rm neigh}(0,\real)$ satisfies (\ref{eq4.1}) for some $k\in \z$. We then claim that, with the norms taken in $L^2$,
\begeq
\label{eq4.4}
\norm{(1-\chi)u}\leq {\cal O}\left(\frac{1}{\eps}\right)\norm{v}+{\cal O}(h^{\infty})\norm{u}.
\endeq
When establishing (\ref{eq4.4}), we only have to prove this bound with $\psi_{1,\pm}$ in place of $1-\chi$, as the estimate involving $\psi_{2,\pm}$ follows
from the semiclassical elliptic regularity.

\bigskip
\noindent
Let $\gamma\subset p^{-1}(0)$ be a closed $H_p$-orbit away from $\widehat{\Lambda}_{0,F_0}$, and assume, to fix the ideas, that
$\langle{q\rangle} \geq F_0 + 1/C$ near $\gamma$. Let $\psi$, $\widetilde{\psi}\in C^{\infty}_0$ be supported in a small flow-invariant neighborhood of
$\gamma$ and assume that $H_p \psi = H_p \widetilde{\psi} = 0$ and that $\widetilde{\psi} = 1$ near ${\rm supp}\, \psi$. In view of a standard iteration
argument~\cite{HS1}, it suffices to prove that
\begeq
\label{eq4.4.1}
\norm{\psi\, u} \leq {\cal O}\left(\frac{1}{\eps}\right)\norm{v}+{\cal O}(h)\norm{\widetilde{\psi}\, u} + {\cal O}(h^{\infty})\norm{u}.
\endeq
In doing so, we shall use the normal form for $P_{\eps}$ near $\gamma$, recalled in the proof of Theorem 1.1 in Section 2. We have
$$
(P_{\eps}-z)\psi\,u = \psi\, v + [P_{\eps},\psi] u.
$$
Here $[P_{\eps},\psi] = {\cal O}(h^3 + \epsilon h) = {\cal O}(\eps h)$ as an operator on $L^2$, since $h^2 \leq \eps$ and the subprincipal symbols of $P_{0}$ and
$\psi$ vanish. Applying the Fourier integral operator $U$ introduced in the proof of Theorem 1.1 and using Egorov's theorem, we obtain, modulo an error term
of norm ${\cal O}(h^{\infty})\norm{u}$,
\begeq
\label{eq4.5}
\left(f(hD_t) + \eps \langle{q\rangle}(hD_t,x,hD_x) + {\cal O}(\eps^2 + h^2)-z\right) U \psi\,u = U\left(\psi\,v + [P_{\eps},\psi] u\right).
\endeq

\medskip
\noindent
Let us now check that the operator $f(hD_t) + \eps \langle{q\rangle}(hD_t,x,hD_x)-z$, acting on $L^2_f(S^1\times \real)$, is invertible, microlocally near
$\tau = x = \xi =0$, with the norm of the inverse being ${\cal O}(1/\eps)$, provided that $z\in {\rm neigh}(0,\real)$ is such that (\ref{eq4.1}) holds. To that
end, we consider a direct sum orthogonal decomposition,
\begin{multline}
\label{eq4.51}
f(hD_t) + \eps \langle{q\rangle}(hD_t,x,hD_x) -z \\
= \bigoplus_{k'\in {\bf Z}} \left(f(h(k'-\theta)) + \eps \langle{q}\rangle(h(k'-\theta),x,hD_x)-z\right),
\end{multline}
where it is understood that we only consider the values of $k'\in \z$ for which $h(k'-\theta)$ is small enough. Using that $\eps \ll h$, we see that
for each $k'\neq k$, with $k$ given in (\ref{eq4.1}), the corresponding direct summand in (\ref{eq4.51}) is invertible, microlocally near $x=\xi =0$ with a
norm of the inverse being ${\cal O}(h^{-1})$. When verifying the microlocal invertibility in the case $k' = k$, we write $z = f(h(k-\theta)) + \eps w$, where
$\abs{w-F_0} \leq 1/C$, $C\gg 1$. We have $\langle{q\rangle}(\tau,x,\xi) - F_0 \sim 1$, for $\tau$, $x$, $\xi \approx 0$, and the operator
$$
f(h(k-\theta)) + \eps \langle{q}\rangle(h(k-\theta),x,hD_x)-z = \eps\left(\langle{q\rangle}(h(k-\theta),x,hD_x)-w\right)
$$
is therefore invertible, microlocally near $x = \xi =0$, with the ${\cal O}(\eps^{-1})$ bound for the norm of the inverse.

\medskip
\noindent
From (\ref{eq4.5}) we therefore infer that
\begeq
\label{eq4.52}
\norm{\psi\, u} \leq {\cal O}\left(\frac{1}{\eps}\right)\left(\norm{v} + \norm{[P_{\eps},\psi] u}\right) + {\cal O}(h^{\infty})\norm{u},
\endeq
and using also that
$$
\norm{[P_{\eps},\psi]u} \leq {\cal O}(\eps h)\norm{\widetilde{\psi}\,u} + {\cal O}(h^{\infty})\norm{u},
$$
we obtain the bounds (\ref{eq4.4.1}) and then (\ref{eq4.4}).

\medskip
\noindent
Relying upon (\ref{eq4.4}), we shall now complete the proof of the fact that the spectrum of $P_{\eps}$ in the region (\ref{eq4.1}) is contained
in the union of the intervals $J_{\ell}(h)$ centered at the quasi-eigenvalues (\ref{eq4.1.1}). Let us write
$$
(P_{\eps}-z)\chi\, u=\chi\, v+[P_{\eps},\chi]u,
$$
where from (\ref{eq4.3}) we know that the norm of the commutator term does not exceed ${\cal O}(h^{\infty})\norm{u}$.
Applying the unitary Fourier integral operator $U$ of Theorem 3.1, we get, modulo an error term of norm ${\cal O}(h^{\infty}) \norm{u}$,
$$
\left(\widehat{P}-z\right)U\chi\, u=U\left(\chi\, v+[P_{\eps},\chi]u\right).
$$
Now an expansion in a Fourier series shows that the operator $\widehat{P}-z$ is invertible, microlocally near $\xi=0$, with a microlocal inverse of the norm
${\cal O}((\eps h)^{-1})$, provided that $z$ in the set (\ref{eq4.1}) avoids the intervals $J_{\ell}(h)$. We get
$$
\norm{\chi u}\leq {\cal O}\left(\frac{1}{\eps h}\right)\norm{v}+{\cal O}(h^{\infty})\norm{u},
$$
and combining this estimate together with (\ref{eq4.4}) we infer that the operator $P_{\eps}-z~:H(m) \rightarrow L^2(M)$
is injective, hence bijective, since it is a Fredholm operator of index zero by general arguments, for $h>0$ small enough.

\bigskip
\noindent
When $z$ in (\ref{eq4.1}) varies in an interval $J_{\ell}(h)$ centered around the quasi-eigenvalue in (\ref{eq4.1.1}), contained in the set in (\ref{eq4.1}),
for some $\ell \in \z$, we may follow Section 6 of~\cite{HS1} and set up a globally well posed Grushin problem for the operator $P_{\eps}-z$. Since the
corresponding discussion here is even simpler than that of~\cite{HS1}, we shall only recall the main steps. Let us define the rank one operators
$$
R_+: L^2(M) \rightarrow \comp,\quad R_- : \comp \rightarrow L^2(M),
$$
given by
$$
R_+ u = (U\chi u, e_{k\ell}),\quad R_- u_- = u_- U^* e_{k\ell}.
$$
Here
$$
e_{k\ell}(x) = \frac{1}{2\pi} e^{\frac{i}{h}(h(k-\theta_1)x_1 + h(\ell-\theta_2)x_2)},\quad \theta_j = \frac{\alpha_j}{4} + \frac{S_j}{2\pi h},\quad j=1,2,
$$
the scalar product in the definition of $R_+$ is taken in the space $L^2_f({\bf T}^2)$, and $U^*$ is the complex adjoint of $U$. The arguments of
Section 6 of~\cite{HS1} can now be applied as they stand to show that for every $(v,v_+)\in L^2(M)\times \comp$, the Grushin problem
$$
(P_{\eps}-z)u + R_- u_- =v,\quad R_+ u = v_+,
$$
has a unique solution $(u,u_-)\in H(m)\times \comp$. We have the corresponding estimate
$$
\eps h\norm{u} + \abs{u_-} \leq {\cal O}(1)\left(\norm{v} + \eps h\abs{v_+}\right),
$$
and if we write the solution in the form
$$
u = Ev + E_+ v_+,\quad u_- = E_- v + E_{-+}v_+,
$$
then repeating the arguments of~\cite{HS1}, we find that
$$
E_{-+}(z)  = z - \widehat{P}\left(h(k-\theta_1), h(\ell -\theta_2),\eps, \frac{h^2}{\eps};h\right) + {\cal O}(h^{\infty}).
$$
Since the eigenvalues of $P_{\eps}$ in the interval $J_{\ell}(h)$ are precisely the values of $z$ for which $E_{-+}(z)$ vanishes~\cite{SjZw}, we see that we
have established Theorem 1.2, in the general case when the clusters of $P_0$ are of size ${\cal O}(h^2)$, and when $\eps$ is in the range $h^2 \ll \eps \ll h$.

\medskip
\noindent
{\it Remark}. The number of the eigenvalues of $P_{\eps}$ in the subcluster (\ref{eq4.1}) is $\sim h^{-1}$, which is of the same order of magnitude as the
total number of eigenvalues of $P_{\eps}$ in the $k$th spectral cluster $f(h(k-\theta)) + [-{\cal O}(\eps), {\cal O}(\eps)]$. See also Chapter 15 of~\cite{DS}.

\section{Improving parameter range for thin clusters}
\setcounter{equation}{0}

In this section, following~\cite{HS3}, we shall extend the range of $\eps$ in Theorem 1.2, in the case when the spectrum of $P_0$ near $0$ clusters into bands of
size ${\cal O}(1)h^{N_0}$, $N_0 > 2$.

\medskip
\noindent
Let $P_{\eps}$, $\eps\in {\rm neigh}(0,\real)$, be a smooth family of selfadjoint operators, such the assumptions of the introduction are satisfied.
As we saw in Section 3, microlocally near the Lagrangian torus $\Lambda_{0,F_0}$, the operator $P_0$ can be reduced by successive averaging procedures to an
operator of the form
\begeq
\label{eq5.1}
P_0 \sim \sum_{j=0}^{\infty} h^j p_j(x_2,\xi),
\endeq
defined near $\xi=0$ in $T^*{\bf T}^2$, and such that $p_0 = p(\xi_1)$, $p_1 =0$. We then have the following result.

\begin{prop}
Assume that the subprincipal symbol of $P_0$ vanishes and that the spectrum of $P_0$ clusters into intervals of size $\leq {\cal O}(h^{N_0})$, for some integer
$N_0 > 2$. Then the terms $p_j(x_2,\xi) = p_j(\xi_1)$ in {\rm (\ref{eq5.1})} are independent of $(x_2,\xi_2)$ when $1\leq j \leq N_0 -1$.
\end{prop}

\medskip
\noindent
Proposition 5.1 is an analog of Proposition 12.1 of~\cite{HS3}, where a microlocal model for the selfadjoint operator $P_0$ near a closed $H_p$--trajectory was
considered. This minor difference does not affect the validity of the result, and the proof of Proposition 5.1 is essentially the same as that of Proposition 12.1
in~\cite{HS3}, making use of a suitable family of ${\cal O}(h^{1/2})$--Gaussian quasimodes on the one-dimensional torus.

\bigskip
\noindent
An application of the discussion in Section 3 together with Proposition 5.1 allows us to conclude that when $0\neq \eps \in {\rm neigh}(0,\real)$, microlocally
near the torus $\Lambda_{0,F_0}$, the operator $P_{\eps}$ can be reduced to the following form,
$$
P_{\eps} = \sum_{j=0}^{\infty} h^j p_j(x_2,\xi,\eps),\quad (x,\xi)\in T^*{\bf T}^2,
$$
where
$$
p_0(x_2,\xi,\eps) = p(\xi_1) + \eps \langle{q\rangle} + {\cal O}(\eps^2)
$$
is independent of $x_1$, and
and
$$
p_1(x_2,\xi,\eps) = \eps q_1(x_2,\xi,\eps),
$$
$$
p_j(x_2,\xi,\eps) = p_j(\xi_1) + \eps q_j(x_2,\xi,\eps),\quad 2 \leq j \leq N_0-1.
$$
It follows that we can write,
\begeq
\label{eq5.2}
P_{\eps} = p(\xi_1;h) + \eps \left(r_0\left(x_2,\xi,\eps,\frac{h^{N_0}}{\eps}\right) + h r_1\left(x_2,\xi,\eps,\frac{h^{N_0}}{\eps}\right) + h^2 r_2\ldots\right),
\endeq
where
$$
p(\xi_1;h) = p(\xi_1) + \sum_{j=2}^{N_0 -1} h^j p_j(\xi_1),
$$
$$
r_0\left(x_2,\xi,\eps,\frac{h^{N_0}}{\eps}\right) = \langle{q\rangle}(\xi) + {\cal O}(\eps) + \frac{h^{N_0}}{\eps}p_{N_0}(x_2,\xi,\eps),
$$
$$
r_1\left(x_2,\xi,\eps,\frac{h^{N_0}}{\eps}\right) = q_1(x_2,\xi,\eps) + \frac{h^{N_0}}{\eps} p_{N_0+1}(x_2,\xi,\eps),
$$
and more generally,
$$
r_j\left(x_2,\xi,\eps,\frac{h^{N_0}}{\eps}\right) = q_j (x_2,\xi,\eps) + \frac{h^{N_0}}{\eps} p_{N_0 + j}(x_2,\xi,\eps),\quad 1 \leq j \leq N_0 -1,
$$
$$
r_j\left(x_2,\xi,\eps,\frac{h^{N_0}}{\eps}\right) = \frac{h^{N_0}}{\eps} p_{j+N_0}(x_2,\xi,\eps),\quad j\geq N_0.
$$

\medskip
\noindent
The analysis of Subsection 3.2 can then be applied to the operator in (\ref{eq5.2}), provided that
$$
\frac{h^{N_0}}{\eps} \leq \delta_0 \ll 1,
$$
and we see that a natural analog of Theorem 3.1 is valid, with the small parameter $h^2/\eps$ replaced by $h^{N_0}/\eps$. The arguments of Section 4 can therefore
also be applied, with minor modifications, and we obtain the full statement of Theorem 1.2, for $\eps$ in the range $h^{N_0}\ll \eps \ll h$.

\section{Magnetic Schr\"odinger operators in the resonant case}
\setcounter{equation}{0}
Let us consider the magnetic Schr\"odinger operator on $\real^2$,
\begeq
\label{eq6.1}
P = \sum_{j=1}^2 (hD_{x_j} + A_j(x))^2 + V(x).
\endeq
Here the magnetic and electric potentials $A = (A_1, A_2)$ and $V$ are assumed to be smooth and real-valued, with
$\partial^{\alpha}A$, $\partial^{\alpha} V\in L^{\infty}(\real^2)$, for all $\alpha\in \nat^2$. It is then well known that $P$ is essentially selfadjoint on
$L^2(\real^2)$, starting from $C_0^{\infty}(\real^2)$.

\bigskip
\noindent
Let us assume that $V\geq 0$ with equality at $0$ only and that $V''(0)>0$. We further assume that
$$
\liminf_{\abs{x}\rightarrow \infty} V(x) > 0.
$$
The spectrum of the selfadjoint nonnegative operator $P$ is then discrete in a neighborhood of $0$.

\medskip
\noindent
Associated to $P$ in (\ref{eq6.1}) is the Weyl symbol given by
\begeq
\label{eq6.2}
p(x,\xi) = \sum_{j=1}^2 (\xi_j + A_j(x))^2 + V(x),\quad x,\xi \in \real^2.
\endeq
Assume that near $0$, for $j=1,2$, we have
\begeq
\label{eq6.3}
A_j(x) = {\cal O}(x^{m-1}),
\endeq
for some $m \geq 3$, and that
\begeq
\label{eq6.4}
V(x) = \frac{1}{2} V''(0)x\cdot x + {\cal O}(x^m).
\endeq
After a linear symplectic change of coordinates, we obtain, as $(x,\xi)\rightarrow 0$,
\begeq
\label{eq6.5}
p(x,\xi) = p_2(x,\xi) + \sum_{j=1}^2 A_{j,m-1}(x) \xi_j + p_m(x) + {\cal O}((x,\xi)^{m+1}).
\endeq
Here
$$
p_2(x,\xi) = \sum_{j=1}^2 \frac{\lambda_j}{2} (x_j^2 + \xi_j^2), \quad \lambda_j > 0,
$$
$A_{j,m-1}$ is a homogeneous polynomial of degree $m-1$, and $p_m(x)$ is a homogeneous polynomial of degree $m$. In what follows, in order to fix the ideas, we shall
consider the case $m = 4$. Assume also, for simplicity, that the electrical potential $V$ satisfies $V(-x) = V(x)$ and the magnetic potential $A$ satisfies
$A(-x) = -A(x)$. We can then rewrite (\ref{eq6.5}) as follows,
\begeq
\label{eq6.5.1}
p(x,\xi) = p_2(x,\xi) + \sum_{j=1}^2 A_{j,3}(x) \xi_j + p_4(x) + {\cal O}((x,\xi)^{6}).
\endeq

\bigskip
\noindent
We assume that $\lambda = (\lambda_1,\lambda_2)$ fulfills the resonant condition,
\begeq
\label{eq6.6}
\lambda \cdot k = 0,
\endeq
for some $0\neq k \in \z^2$. We shall then be interested in eigenvalues $E$ of $P$ with $E\sim \eps$, where $h^{2\delta} < \eps \ll 1$, $0 < \delta < 1/2$. The
general arguments of~\cite{Sj92} imply that the corresponding eigenfunctions are microlocally concentrated in the region where $(x,\xi) = {\cal O}(\eps^{1/2})$,
and we introduce therefore the change of variables $x = \eps^{1/2}y$. Then
$$
\frac{1}{\eps} P(x,hD_x) = \frac{1}{\eps}P(\eps^{1/2} (y,\widetilde{h}D_y)),\quad \widetilde{h} = \frac{h}{\eps} \ll 1.
$$
It follows from (\ref{eq6.5.1}) that the symbol of the corresponding $\widetilde{h}$-pseudodifferential operator is
$$
\frac{1}{\eps} p(\eps^{1/2}(y,\eta)) = p_2(y,\eta) + \eps q(y,\eta) + {\cal O}(\eps^{2}),
$$
to be considered in the region where $\abs{(y,\eta)} = {\cal O}(1)$. Here
\begeq
\label{eq6.7}
q(y,\eta) = \sum_{j=1}^2 A_{j,3}(y) \eta_j + p_4(y).
\endeq
The resonant assumption (\ref{eq6.6}) implies that the $H_{p_2}$--flow is periodic on $p_2^{-1}(E)$, for $E\in {\rm neigh}(1,\real)$, with period $T>0$ which
does not depend on $E$, and we shall assume that $T$ is the minimal period for the $H_{p_2}$--flow. We may therefore apply Theorem 1.2 to discuss the
invertibility of
$$
P(x,hD_x) - \eps(1+z) = \eps\left(\frac{1}{\eps} P(x,hD_x) -1-z\right) ,\quad z\in {\rm neigh}(0,\real),
$$
in the range of energies $E = \eps(1+z)$, given by
$$
h^{N_0/(N_0 + 1)} \ll E \ll h^{1/2},
$$
for all $N_0 = 2,3,\ldots\,$. Notice also that since the eigenvalues of $p^w_2(x,\widetilde{h}D_x)$ depend linearly on $\widetilde{h}$, the functions $f_j$,
$j\geq 2$, occurring in Theorem 1.2, all vanish. We obtain the following result.

\begin{prop}
Assume that {\rm (\ref{eq6.6})} holds and that the $H_{p_2}$--flow has a minimal period $T>0$ on $p_2^{-1}(1)$. Let $\langle{q\rangle}$ stand for the average of
the homogeneous function $q$ in {\rm (\ref{eq6.7})} along the trajectories of the Hamilton vector field of $p_2$, and assume that $\langle{q\rangle}$ is not
identically zero. Let $F_0\in \real$ be a regular value of $\langle{q\rangle}$ restricted to $p^{-1}_2(1)$. Let $\eps$ satisfy
$$
h^{N_0/N_0 + 1}\ll \eps \ll h^{1/2},
$$
for some $N_0 \geq 2$ fixed. Then for $z\in {\rm neigh}(0,\real)$ in the set
$$
\abs{z - f\left(\widetilde{h}(k - \frac{\alpha_1}{4})-\frac{S_1}{2\pi}\right) - \eps F_0} < \frac{\eps}{{\cal O}(1)},\quad f(E) = \frac{2\pi}{T} E, \quad
\widetilde{h} = \frac{h}{\eps},
$$
the eigenvalues of $P$ of the form $\eps(1+z)$ are given by
$$
z=\widehat{P}\left(\widetilde{h}(k-\frac{\alpha_1}{4})-\frac{S_1}{2\pi},\widetilde{h}(\ell-\frac{\alpha_2}{4})-\frac{S_2}{2\pi},\eps,
\frac{\widetilde{h}^{N_0}}{\eps};\widetilde{h}\right) + {\cal O}(h^{\infty}),\quad \ell \in \z.
$$
Here $\widehat{P}(\xi,\eps,\widetilde{h}^{N_0}/\eps;\widetilde{h})$ has an expansion, as $\widetilde{h}\rightarrow 0$,
$$
\widehat{P}\left(\xi,\eps,\frac{\widetilde{h}^{N_0}}{\eps};\widetilde{h}\right)\sim f(\xi_1)+\eps \sum_{n=0}^{\infty} \widetilde{h}^n
r_n\left(\xi,\eps,\frac{\widetilde{h}^{N_0}}{\eps}\right),
$$
where
$$
r_0(\xi)=\langle{q}\rangle(\xi)+{\cal O}\left(\eps+\frac{\widetilde{h}^{N_0}}{\eps}\right),\quad r_j={\cal O}(1),\quad j\geq 1.
$$
The coordinates $\xi_1 = \xi_1(E)$ and $\xi_2 = \xi_2(E,F)$ are the normalized actions of the Lagrangian tori
$$
\Lambda_{E,F}: p_2 = E, \quad \langle{q\rangle} = F,
$$
for $E\in {\rm neigh}(1,\real)$, $F\in {\rm neigh}(F_0,\real)$, given by
$$
\xi_j = \frac{1}{2\pi} \left(\int_{\gamma_j(E,F)}\eta\,dy - \int_{\gamma_j(1,F_0)} \eta\,dy\right),\quad j = 1,2,
$$
with $\gamma_j(E,F)$ being fundamental cycles in $\Lambda_{E,F}$, such that $\gamma_1(E,F)$ corresponds to a closed $H_{p_2}$--trajectory of minimal period $T$.
Furthermore,
$$
S_j = \int_{\gamma_j(1,F_0)} \eta\,dy,
$$
and $\alpha_j\in \z$ is fixed, $j=1,2$.
\end{prop}

\medskip
\noindent
We shall finish this section by providing an explicit example, illustrating Proposition 6.1 in the case when $\lambda = (1,1)$. Then $T= 2\pi$ is the minimal
period for the $H_{p_2}$--flow, and our task becomes computing the flow average $\langle{q\rangle}$ and determining its critical values, viewed as a function on the
compact symplectic manifold $\Sigma$. In this case, as we saw in~\cite{HS3}, the manifold $\Sigma$ can naturally be identified with the complex projective line
$\comp P^1 \cong S^2$.

\medskip
\noindent
Continuing to follow~\cite{HS3}, let us recall first how to compute the trajectory average of a monomial $x^\alpha \xi ^\beta $ with $\abs{\alpha} +\abs{\beta}=m$,
for some $m\in\{3,4,5,...\}$. To this end, it is convenient to introduce
\begeq
\label{eq6.8}
z_j=x_j+i\xi _j\in \comp,\quad j=1,2,
\endeq
and we then notice that along a $H_{p_2}$-trajectory we get in the $z_1,z_2$ coordinates:
\begeq
\label{eq6.9}
z_j(t)=e^{-i\lambda_jt}z_j(0).
\endeq
Then we write $x_j(t)=\Re z_j(t)$, $\xi _j(t)=\Im z_j(t)$, so that
\begin{multline}
\label{eq6.81}
x(t)^\alpha \xi (t)^\beta = \prod_{j=1}^2 ((\Re z_j(t))^{\alpha _j}(\Im
z_j(t))^{\beta _j}) \\
={1\over 2^{\vert \alpha \vert +\vert \beta \vert }i^{\vert \beta \vert }
}\prod_{j=1}^2 ((z_j(0)e^{-i\lambda _jt}+\overline{z_j}(0) e^{i\lambda
_jt})^{\alpha _j}(z_j(0)e^{-i\lambda _jt}-\overline{z_j}(0) e^{i\lambda
_jt})^{\beta _j}).
\end{multline}
Expanding the product by means of the binomial theorem, we see that the time average is equal to the time-independent
term, and since this average is constant along each trajectory we shall replace the symbols $z_j(0)$ simply by $z_j$.

\medskip
\noindent
For simplicity, we shall assume that $p_4 = 0$ in (\ref{eq6.7}), and then we write
\begeq
\label{eq6.10}
A_{j,3}(x) = \sum_{k=0}^3 a_{j,k} x_1^k x_2^{3-k}, \quad j =1,2.
\endeq
The associated magnetic field
$$
B(x) = \frac{\partial A_{2,3}}{\partial{x_1}} - \frac{\partial A_{1,3}}{\partial{x_2}}
$$
is given by
\begeq
\label{eq6.11}
B(x) = b_2 x_1^2 + b_1 x_1 x_2 + b_0 x_2^2,
\endeq
where
\begeq
\label{eq6.12}
b_2 = 3 a_{2,3} - a_{1,2},\quad b_1 = 2(a_{2,2} - a_{1,1}), \quad b_0 = a_{2,1} - 3a_{1,0}.
\endeq

\medskip
\noindent
Using (\ref{eq6.81}), we get
$$
\langle{\xi_1 x_1^3\rangle} = 0,
$$
$$
\langle{\xi_1 x_1^2 x_2\rangle} = \frac{1}{16 i} (z_1^2 \overline{z}_1\overline{z}_2 - \overline{z}_1^2 z_1 z_2) = -\frac{1}{2} \rho_1^{3/2} \rho_2^{1/2}
\sin(\theta_1 - \theta_2),
$$
$$
\langle{\xi_1 x_1 x_2^2\rangle} = \frac{1}{16 i} (z_1^2 \overline{z}_2^2 - \overline{z}_1^2 z_2^2) = -\frac{1}{2} \rho_1 \rho_2
\sin 2(\theta_1 - \theta_2),
$$
$$
\langle{\xi_1 x_2^3\rangle} = \frac{3}{16 i} (z_1 z_2 \overline{z}_2^2 - \overline{z}_1 \overline{z}_2 z_2^2) = -\frac{3}{2} \rho_1^{1/2} \rho_2^{3/2}
\sin (\theta_1 - \theta_2),
$$
$$
\langle{\xi_2 x_1^3\rangle} = \frac{3}{16 i} (z_1 z_2 \overline{z}_1^2 - \overline{z}_1 \overline{z}_2 z_1^2) = \frac{3}{2} \rho_1^{3/2} \rho_2^{1/2}
\sin (\theta_1 - \theta_2),
$$
$$
\langle{\xi_2 x_1^2 x_2\rangle} = \frac{1}{16 i} (z_2^2 \overline{z}_1^2 - \overline{z}_2^2 z_1^2) = \frac{1}{2} \rho_1 \rho_2
\sin2(\theta_1 - \theta_2),
$$
$$
\langle{\xi_2 x_1 x_2^2\rangle} = \frac{1}{16 i} (z_2^2 \overline{z}_1\overline{z}_2 - \overline{z}_2^2 z_1 z_2) = \frac{1}{2} \rho_1^{1/2} \rho_2^{3/2}
\sin(\theta_1 - \theta_2),
$$
$$
\langle{\xi_2 x_2^3\rangle} = 0.
$$
Here $(\rho_j,\theta_j)$ are the action--angle variables given by
$$
z_j = \sqrt{2\rho_j} e^{-i\theta_j}.
$$
Recalling the expressions for $q$ in (\ref{eq6.7}) and $B$ in (\ref{eq6.11}), (\ref{eq6.12}), we get
\begeq
\label{eq6.13}
2\langle{q\rangle} = b_2 \rho_1^{3/2} \rho_2^{1/2} \sin(\theta_1 - \theta_2) + \frac{b_1}{2} \rho_1 \rho_2 \sin2(\theta_1 - \theta_2) +
b_0 \rho_1^{1/2} \rho_2^{3/2} \sin(\theta_1 - \theta_2).
\endeq
In particular, as was already observed in the example in the introduction, the flow average $\langle{q\rangle}$ depends on the magnetic field $B$ only. Let us
consider the special case when $b_2 = b_0 = 0$ while $b_1 \neq 0$. In this case, a straightforward computation shows that $\langle{q\rangle}$, viewed as a function
on the space $\Sigma$, has exactly three critical values, given by $\pm b_1/16$ and $0$. When $F_0 \in \real$ is in the range of $\langle{q\rangle}$, $F_0$ away
from $\pm b_1/8$ and $0$, Proposition 6.1 applies.

\end{document}